\documentclass [11 pt,a4paper,twoside]{article}

\usepackage[english]{babel}
\usepackage[T1]{fontenc}
\usepackage[latin1]{inputenc}

\usepackage{geometry}

\usepackage{enumerate}

\usepackage{graphicx}
\usepackage{pst-tree}
\usepackage{pstricks}
\usepackage{pst-eucl}

\usepackage{amsmath,amsfonts,amssymb,mathrsfs}
\usepackage{theorem}
\usepackage{textcomp}
\usepackage{float}

\newtheorem{theo}{Theorem}
\newtheorem{corol}{Corollary}
\newtheorem{prop}{Proposition}
\newtheorem{lemm}{Lemma}

\newcommand{\pen}{\text{\upshape{pen}}}
\newcommand{\argmin}[1]{\underset{#1}{\text{argmin }}}
\newcommand{\argmintext}[1]{{\text{argmin}}_{#1}}
\newcommand{\Var}{\text{\upshape{Var}}}

\newcommand{\bb}[1]{\boldsymbol{#1}}

\newdimen\AAdi%
\newbox\AAbo%
\def\AArm{\fam0 }
\def\AAk#1#2{\setbox\AAbo=\hbox{#2}\AAdi=\wd\AAbo\kern#1\AAdi{}}%
\def\BBone{{\AArm 1\AAk{-.8}{I}I}}%

\renewcommand{\d}{\mathrm{d}}
\newcommand{\N}{\mathbb{N}}

\newcommand{\R}{\mathbb{R}}
\renewcommand{\P}{\mathbb{P}}
\newcommand{\E}{\mathbb{E}}
\renewcommand{\L}{\mathbb{L}}
\newcommand{\mr}{{(m,\bb r)}}
\newcommand{\mrho}{{(m,\bb \rho)}}
\newcommand{\primemrho}{{(m',\bb \rho')}}
\newcommand{\hatmrho}{{(\hat m,\bb{\hat\rho})}}
\newcommand{\starmr}{{(m_\star,\bb{r_\star})}}
\newcommand{\M}{\mathcal M_\star}
\newcommand{\Mdeg}{\boldsymbol{\mathcal M_\star^{deg}}}
\newcommand{\Mone}{{M_{1,\star}}}

\title{Adaptation to anisotropy and inhomogeneity via dyadic piecewise polynomial selection}

\author{Nathalie {\scshape Akakpo}  \medskip\\ \textit{LPMA, Universit\'e Pierre et Marie Curie}}

\begin{document}
\maketitle

\begin{abstract}
This article is devoted to nonlinear approximation and estimation via piecewise polynomials built on partitions into dyadic rectangles. The approximation rate is studied over possibly inhomogeneous and anisotropic smoothness classes that contain Besov classes. Highlighting the interest of such a result in statistics, adaptation in the minimax sense to both inhomogeneity and anisotropy of a related multivariate density estimator is proved. Besides, that estimation procedure can be implemented with a computational complexity simply linear in the sample size.
\end{abstract}



\tableofcontents

\section{Introduction}

When estimating a multivariate function, it seems natural to consider that its smoothness is likely to vary either spatially, or with the direction, or both. We will refer to the first feature as (spatial) inhomogeneity. If the risk is measured in a $\L_q$-norm, measuring the smoothness in a $\L_p$-norm with $p<q$ allows to take into account such an inhomogeneity -- all the greater as $p$ is smaller -- in the sense that functions with some localized singularities and otherwise flat parts may thus keep a high smoothness index. For the second feature, we will talk about anisotropy, which is usually described by different indices of smoothness according to the coordinate directions. Yet, statistical procedures that adapt both to possible inhomogeneity and anisotropy remain rather scarce. Indeed, the existing literature seems to amount to the following references. Neumann and Von Sachs~\cite{NeumannVonSachs}, for estimating the evolutionary spectrum of a locally stationary time series, and Neumann~\cite{Neumann}, in the Gaussian white noise framework, study thresholding procedures in a tensor product wavelet basis. In a Gaussian regression framework, Donoho~\cite{DonohoCart} proposes the dyadic CART procedure, a selection procedure among histograms built on partitions into dyadic rectangles, extended to the density estimation framework by Klemel\"a~\cite{Klemela}. Last, Kerkyacharian, Lepski and Picard~\cite{KLP} introduce a kernel estimator with adaptive bandwidth in the Gaussian white noise model. These authors study the performance of their procedures for the $\L_2$-risk, apart from the latter who consider any $\L_q$-risk for $q\geq 1$. Neumann and Von Sachs~\cite{NeumannVonSachs} measure the smoothness of the function to estimate in the Sobolev scale, whereas the others consider the finer Besov scale. Besides, the $\L_p$-norm in which the smoothness is measured is allowed to vary with the direction, except in~\cite{DonohoCart}, but always constrained to be greater than 1. Common to those few procedures is the ability to reach the minimax rate over a wide range of possibly inhomogeneous and anisotropic classes, up to a logarithmic factor, the unknown smoothness being as usually limited by the \textit{a priori} fixed smoothness of the underlying wavelets, piecewise polynomials or kernel.

Adaptation results of the aforementioned type rely as much on Statistics as on Approximation Theory, oracle-type inequalities reflecting the interplay between both domains. Assume for instance that the function $s$ to estimate lies in the set $\mathcal F([0,1]^d,\R)$ of all real-valued functions defined over the unit cube $[0,1]^d$, let $(S_m)_{m\in\mathcal M}$ be a given family of linear subspaces of $\mathcal F([0,1]^d,\R)$ and $\tilde s$ a statistical procedure somehow based on that family. An oracle-type inequality in the $\L_q$-norm roughly takes the form
\begin{equation}
\E_s\left[\|s-\tilde s\|_q^q\right] \leq C \inf_{m\in\mathcal M} \left\{\inf_{t\in S_m} \|s-t\|_q^q + \left(\dim(S_m)/n\right)^{q/2}\right\},
\label{eq:oracletype}
\end{equation}
where $C$ is some positive constant, indicating that $\tilde s$ is able to choose a model $S_m$ in the family that approximately realizes the best compromise between the approximation error and the dimension of the model. Equivalently, it may be written as 
\begin{equation}
\E_s\left[\|s-\tilde s\|_q^q\right] \leq C \inf_{D\in\N^\star}\left\{\inf_{t\in\cup_{m\in\mathcal M_D}S_m} \|s-t\|_q^q + \left(D/n\right)^{q/2}\right\},
\label{eq:oracletype2}
\end{equation}
where $\mathcal M_D=\{m\in\mathcal M \text{ s.t. } \dim(S_m)=D\}.$ On the other hand, the collection $(S_m)_{m\in\mathcal M}$ should be chosen so as to have good approximation properties over various classes $\mathcal S(\alpha,p,R)$ of functions with smoothness $\alpha$ measured in a $\L_p$-norm and with semi-norm smaller than $R$. Otherwise said, each approximating space $\cup_{m\in\mathcal M_D}S_m$ -- typically nonlinear to deal with inhomogeneous functions-- should satisfy, for a wide range of values of $\alpha, p$ and $R$,
\begin{equation}
\sup_{s\in\mathcal S(\alpha,p,R)} \inf_{t\in\cup_{m\in\mathcal M_D}S_m} \|s-t\|_q \leq  C(\alpha,p) RD^{- \alpha/d},
\label{eq:adaptiveapprox}
\end{equation}
for some positive real $C(\alpha,p)$ that only depends on $\alpha$ and $p$. Combining the oracle-type inequality~\eqref{eq:oracletype2} and the approximation result~\eqref{eq:adaptiveapprox} then provides an estimator $\tilde s$ with rate at most of order $(Rn^{-\alpha/d})^{qd/(d+2\alpha)}$ over each class $\mathcal S(\alpha,p,R)$, which is usually the minimax rate. Having at one's disposal spaces $(S_m)_{m\in\mathcal M}$ that do no depend on any \textit{a priori} knowledge about the smoothness of the function to estimate -- other than the scale of spaces it belongs to -- and reaching the approximation rate~\eqref{eq:adaptiveapprox} is thus a real issue for statisticians. In order to deal with inhomogeneity only, in a multivariate framework, such results appear for instance in the following references. DeVore, Jawerth and Popov~\cite{DeVoreJP}, Birg\'e and Massart~\cite{BMBesov} or Cohen, Dahmen, Daubechies and  DeVore~\cite{CohenTree}  propose wavelet based approximation algorithms aimed in particular at Besov type smoothness. Applications of the approximation result of~\cite{BMBesov}  to statistical estimation may be found in Birg\'e and Massart~\cite{BM97} or Massart~\cite{Massart} for instance. DeVore and Yu~\cite{DeVore} are concerned with piecewise polynomials built on partitions into dyadic cubes, notably for functions with Besov type smoothness. But their result will wait until Birg\'e~\cite{BirgeBrouwer} to be used in Statistics. More generally, such results are in fact hidden behind all adaptive procedures. Thus, for both inhomogeneous and anisotropic functions, we refer in particular to the articles cited in the first paragraph. Let us underline that the procedure studied by Donoho~\cite{DonohoCart} and Klemel\"a~\cite{Klemela}, though based on dyadic rectangles instead of cubes, does not rely on a nonlinear approximation result via piecewise polynomials such as~\cite{DeVoreYu}. Indeed, the adaptivity of that estimator follows from its characterization as a wavelet selection procedure among some tree-structured subfamily of the Haar basis. Other nonlinear wavelet based approximation results are proved in Hochmuth~\cite{HochmuthNonlinear} or~\cite{Leisner} for anisotropic Besov spaces. Last, piecewise constant approximation based on dyadic rectangles is studied in Cohen and Mirebeau~\cite{CohenMirebeau} for nonstandard smoothness spaces under the constraint of continuous differentiability. 

Our aim here is to provide an approximation result tailored for statisticians, whose interest is illustrated by a new statistical procedure. The first part of the article is devoted to piecewise polynomial approximation based on partitions into dyadic rectangles. Thanks to an approximation algorithm inspired from DeVore and Yu~\cite{DeVoreYu}, we obtain approximation rates akin to~\eqref{eq:adaptiveapprox} over possibly inhomogeneous and anisotropic smoothness classes that contain for instance the more traditional Besov classes. The approximation rate can be measured in any $\L_q$-norm, for $1\leq q\leq \infty$, and we allow an arbitrarily high inhomogeneity in the sense that we measure the smoothness in a $\L_p$-norm with $p$ allowed to be arbitrarily close to 0. Besides, we take into account a possible restriction on the minimal size of the dyadic rectangles, which may arise in statistical applications. For estimating a multivariate function, we then introduce a selection procedure that chooses from the data the best partition into dyadic rectangles and the best piecewise polynomial built on that partition thanks to a penalized least-squares type criterion. The degree of the polynomial may vary from one rectangle to another, and also according to the coordinate directions, so as to provide a good adaptation both to inhomogeneity and anisotropy. Thus, our procedure extends the dyadic histogram selection procedures of Donoho~\cite{DonohoCart}, Klemel\"a~\cite{Klemela} or Blanchard, Sch\"afer, Rozenholc and M\"uller~\cite{BSRM}, and the dyadic piecewise polynomial estimation procedure proposed in a univariate or isotropic framework by Willett and Nowak~\cite{WillettNowak}. We study the theoretical performance of the procedure -- with no need to resort to the "wavelet trick" used in~\cite{DonohoCart,Klemela} -- for the $\L_2$-risk in the density estimation framework, as~\cite{Klemela}, but we propose a more refined form of penalty than~\cite{Klemela}. For such a penalty, we provide an oracle-type inequality and adaptivity results in the minimax sense over a wide range of possibly inhomogeneous and anisotropic smoothness classes that contain Besov type classes. We emphasize that, if the maximal degree of the polynomials does not depend on the sample size, we reach the minimax rate up to a constant factor only, contrary to \textit{all} the previously mentioned estimators. This results not only from the good approximation properties of dyadic piecewise polynomials, but also from the moderate number of dyadic partitions of the same size. We can also allow the maximal degree of the polynomials to grow logarithmically with the sample size, in which case we reach the minimax rate on a growing range of smoothness classes, up to a logarithmic factor. Moreover, our procedure can be implemented with a computational complexity only linear in the sample size, possibly up to a logarithmic factor, depending on the way we choose the maximal degree.

The plan of the paper is as follows. Section~\ref{sec:approximation} is devoted to piecewise polynomial approximation based on partitions into dyadic rectangles. In Section~\ref{sec:estimation}, we are concerned with density estimation based on a data-driven choice of a best dyadic piecewise polynomial. We study there the theoretical properties of the procedure and briefly describe the algorithm to implement it. Most proofs of Sections~\ref{sec:approximation}  and~\ref{sec:estimation}  are deferred respectively to Section~\ref{sec:approxproof} and to Sections~\ref{sec:proofselect} and~\ref{sec:proofadapt}.

\section{Adaptive approximation by dyadic piecewise polynomials}\label{sec:approximation}

In this section, we present an approximation algorithm by piecewise polynomials built on partitions into dyadic rectangles. We study its rate of approximation over some classes of functions that may present at the same time anisotropic and inhomogeneous smoothness.

\subsection{Notation}
Throughout the article, we fix $d\in\N^\star$, and throughout this section, we fix some $d$-uple of nonnegative integers $\bb r=(r_1,\ldots,r_d)$ that represent the maximal degree of polynomial approximation in each direction. We call dyadic rectangle of $[0,1]^d$ any set of the form $I_1\times\ldots\times I_d$ where, for all $1\leq l\leq d$,  
$$I_l=[0,2^{-{j_l}}] \quad\text{ or }\quad
I_l=]k_l 2^{-{j_l}},(k_l+1)2^{-{j_l}}]$$
with $j_l \in \N$ and $k_l\in\{1,\ldots,2^{j_l}-1\}$. Otherwise said, a dyadic rectangle of $[0,1]^d$ is defined as a product of $d$ dyadic intervals of $[0,1]$ that may have different lengths. For a partition $m$ of $[0,1]^d$ into dyadic rectangles, we denote by $|m|$ the number of rectangles in $m$ and by $S_{(m,\bb r)}$ the space of all piecewise polynomial functions on $[0,1]^d$ which are polynomial with degree $\leq r_l$ in the $l$-th direction, $l=1,\ldots,d$, over each rectangle of $m$. Besides, for $0< p \leq \infty$, we denote by $\L_p([0,1]^d)$ the set of all real-valued and measurable functions $s$ on $[0,1]^d$ such that the (quasi-)norm
\begin{equation*}
\|s\|_p= 
\begin{cases}
\left(\int_{[0,1]^d} |s(x)|^p \d \lambda_d (x)\right)^{1/p} &\text{ if } 0<p<\infty\\
\sup_{x\in [0,1]^d} |s(x)| &\text{ if } p=\infty
\end{cases}
\end{equation*}
is finite, where $\lambda_d$ is the Lebesgue measure on $[0,1]^d$. Last, $C(\theta)$, $C_i(\theta)$ or $C'_i(\theta)$, $i\in\N^\star$ stand for a positive reals that only depend on the parameter $\theta$. Their values may change from one line to another, unless otherwise said.

\subsection{Approximation algorithm}\label{sec:approxalgo}
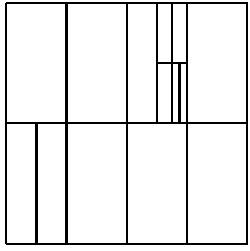
\begin{figure}[h]
\setlength{\unitlength}{0.1cm}
\begin{center}
\begin{picture}(32,32)
\put(0,0){\line(1,0){32}}
\put(0,0){\line(0,1){32}}
\put(0,32){\line(1,0){32}}
\put(32,32){\line(0,-1){32}}
\put(8,0){\line(0,1){32}}
\put(16,0){\line(0,1){32}}
\put(24,0){\line(0,1){32}}
\put(4,0){\line(0,1){16}}
\put(20,16){\line(0,1){16}}
\put(22,16){\line(0,1){16}}
\put(23,16){\line(0,1){8}}
\put(0,16){\line(1,0){32}}
\put(20,24){\line(1,0){4}}
\end{picture}
\caption{Example of partition of $[0,1]^2$ into dyadic rectangles from $\mathcal D^{\bb \sigma}$ for $\sigma_2=2\sigma_1$.}	
\label{flottants:dyadpartapprox:figure}        
\end{center}
\end{figure}

Let us fix $1\leq q \leq \infty$. In order to approximate a possibly anisotropic and inhomogeneous function $s$ in the $\L_q$-norm, we propose an approximation algorithm inspired from~\cite{DeVoreYu}. We shall construct an adequate piecewise polynomial approximation on a partition into dyadic rectangles adapted to $s$, beginning with the trivial partition of the unit square $[0,1]^d$ and proceeding to successive refinements. For doing so, we consider the criterion  
\begin{equation}
\mathcal{E}_{\bb r} (s,K)_q=\inf_{P\in\mathscr{P}_{\bb r}} \|(s-P)\BBone_K\|_q
\label{eq:criterion}
\end{equation} 
measuring the error in approximating $s$ on a rectangle $K\subset [0,1]^d$ by some element from the set $\mathscr{P}_{\bb r}$ of all polynomials on $[0,1]^d$ with degree $\leq r_l$ in the $l$-th direction. We also fix some threshold $\epsilon>0$ -- to be chosen later, according to the smoothness assumptions on $s$. But contrary to~\cite{DeVoreYu}, we allow the degrees of smoothness of $s$ to vary with the directions and describe them by a multi-index $\bb \sigma=(\sigma_1,\ldots,\sigma_d) \in \prod_{l=1}^d (0,r_l+1)$, in a sense that will be made precise in the next subsection. Thus, our algorithm is based on a special subcollection of dyadic rectangles adapted to an anisotropic smoothness measured by $\bb \sigma$. Indeed, for $j\in\N$, we define $\mathcal{D}^{\bb \sigma}_j$ as the set of all dyadic rectangles $I_1\times\ldots\times I_d\subset[0,1]^d$ such that, for all $1\leq l\leq d$, 
$$I_l=\left[0,2^{-\lfloor j\underline{\bb \sigma}/\sigma_l\rfloor}\right]
\quad\text{ or }\quad
I_l=\left]k_l2^{-\lfloor j \underline{\bb \sigma}/\sigma_l \rfloor},
(k_l+1)2^{-\lfloor j\underline{\bb \sigma}/\sigma_l\rfloor}\right],$$
with $\underline{\bb \sigma} =\min_{1\leq l\leq d} \sigma_l$ and $k_l\in\{1,\ldots, 2^{\lfloor j\underline{\bb \sigma}/\sigma_l \rfloor}-1\}$, and we set 
$\mathcal{D}^{\bb \sigma}=\cup_{j\in\N} \mathcal{D}_j^{\bb \sigma}.$
It should be noticed that, for all $j\in\N$, any $K\in\mathcal{D}^{\bb \sigma}_j$ can be partitioned into dyadic rectangles of $\mathcal{D}^{\bb \sigma}_{j+1}$, that we call children of $K$. For $d=2$ and $\sigma_2=2\sigma_1$ for instance, a partition of $[0,1]^2$ into dyadic rectangles from $\mathcal{D}^{\bb \sigma}$ will thus be roughly twice as fine in the first direction, as illustrated by Figure~\ref{flottants:dyadpartapprox:figure}.

The algorithm begins with the set $\mathcal{I}^1(s,\epsilon)$ that only contains $[0,1]^d$. If $\mathcal{E}_{\bb r}(s,[0,1]^d)_q<\epsilon$, then the algorithm stops. Else, $[0,1]^d$ is replaced with his children in $\mathcal{I}^1(s,\epsilon)$, hence a new partition $\mathcal{I}^2(s,\epsilon)$. In the same way, the $k$-th step begins with a partition $\mathcal{I}^k(s,\epsilon)$ of $[0,1]^d$ into dyadic rectangles that belong to $\mathcal{D}^{\bb \sigma}$. If $\max_{K\in\mathcal{I}^k(s,\epsilon)} \mathcal{E}_{\bb r}(s,K)_q <\epsilon,$ then the algorithm stops. Else, a dyadic rectangle $K\in\mathcal{I}^k(s,\epsilon)$ such that $\mathcal{E}_{\bb r}(s,K)_q \geq \epsilon$ is chosen and replaced with his children in $\mathcal{I}^k(s,\epsilon)$, hence a new partition $\mathcal{I}^{k+1}(s,\epsilon)$. Since $s\in\L_q([0,1]^d)$, $\mathcal{E}_{\bb r}(s,K)_q$ tends to 0 when the Lebesgue measure of $K$ tends to 0, so the algorithm finally stops. The final partition $\mathcal{I}(s,\epsilon)$ only contains dyadic rectangles that belong to $\mathcal{D}^{\bb \sigma}$ and such that $\max_{K\in\mathcal{I}(s,\epsilon)} \mathcal{E}_{\bb r}(s,K)_q <\epsilon.$ 
For all $K\in\mathcal{I}(s,\epsilon)$, we approximate $s$ on $K$ by $Q_K(s)$, a polynomial function with degree $\leq r_l$ in the $l$-th direction such that $\|(s-Q_K(s))\BBone_K\|_q=\mathcal{E}_{\bb r}(s,K)_q.$ Otherwise said, we approximate $s$ on the unit cube by   
$$A(s,\epsilon)=\sum_{K\in\mathcal{I}(s,\epsilon)}Q_K(s),$$
thus committing the error  
\begin{equation}
\|s-A(s,\epsilon)\|_q
=\left(\sum_{K\in\mathcal{I}(s,\epsilon)} \|(s -Q_K(s))\BBone_K\|_q^q\right)^{1/q}
< |\mathcal{I}(s,\epsilon)|^{1/q} \epsilon 
\label{eq:simpleapproxfinite}
\end{equation}
if $1\leq q <\infty$, and 
\begin{equation}
\|s-A(s,\epsilon)\|_\infty
=\max_{K\in\mathcal{I}(s,\epsilon)} \|(s -Q_K(s))\BBone_K\|_\infty
< \epsilon
\label{eq:simpleapproxinfinite}
\end{equation}
if $q=\infty.$

\subsection{Approximation rate over anisotropic function classes}\label{sec:approxrate}
In order to study the approximation rate of the previous algorithm, we introduce function spaces that arise naturally from the way the algorithm proceeds. Let us fix $\bb \sigma\in\prod_{l=1}^d (0,r_l+1)$ and $0<p,p'\leq\infty$. For $s\in\L_p([0,1]^d)$ and $k\in\N$, we set
\begin{equation}
e_{\bb r,\bb \sigma,p,k}(s)=\inf_{P\in\Pi^{ \bb r,\bb \sigma}_k} \|s-P\|_p
\label{eq:pperrork}
\end{equation}
where $\Pi^{\bb r,\bb \sigma}_k$ is the set of all piecewise polynomial functions on $[0,1]^d$ that are polynomial with degree $\leq r_l$ in the $l$-th direction over each rectangle in $\mathcal{D}_k^{\bb \sigma}$. Then, we define $\mathcal N^{\bb r,\bb \sigma}_{p'}(\L_p([0,1]^d))$ as the set of all functions $s\in\L_p([0,1]^d)$ such that the quantity 
\begin{equation*}
N_{\bb r,\bb \sigma,p,p'}(s)=
\begin{cases}
\left(\sum_{k\in\N} \left(2^{k\underline{\bb \sigma}} e_{\bb r,\bb \sigma,p,k}(s)\right)^{p'}\right)^{1/p'}& \text{ if } 0<p'<\infty\\
\sup_{k\in\N} \left(2^{k\underline{\bb \sigma}} e_{\bb r,\bb \sigma,p,k}(s)\right) & \text{ if } p'=\infty
\end{cases}
\end{equation*}
is finite. One can easily verify that $N_{\bb r,\bb \sigma,p,p'}$ is a (quasi-)semi-norm on $\mathcal N^{\bb r,\bb \sigma}_{p'}(\L_p([0,1]^d))$, and that $\mathcal N^{\bb r,\bb \sigma}_{p'}(\L_p([0,1]^d))$ gets larger as $p'$ increases since 
\begin{equation}
N_{\bb r,\bb \sigma,p,p_2'}(s) \leq N_{\bb r,\bb \sigma,p,p_1'}(s) \text{ for } 0< p'_1 \leq p'_2\leq \infty.
\label{eq:embeddingNN}
\end{equation}
If $p\geq q$, then $\mathcal N^{\bb r,\bb \sigma}_{p'}(\L_p([0,1]^d))$ is obviously embedded in the space $\L_q([0,1]^d)$ in which we measure the quality of approximation. The same property still holds for $p$ smaller than $q$, under adequate assumptions on the harmonic mean $H(\bb \sigma)$ of $\sigma_1,\ldots,\sigma_d$, \textit{i.e.}  
$$H(\bb \sigma)=\left(\frac{1}{d} \sum_{l=1}^d\frac{1}{\sigma_l}\right)^{-1}.$$
Indeed, denoting by $(x)_+=\max\{x,0\}$ for any real $x$, we prove in Section~\ref{sec:approxproof} the following continuous embedding.
\begin{prop}\label{embeddingLq}
Let $\bb \sigma\in\prod_{l=1}^d (0,r_l+1)$, $0<p,p'\leq \infty$ and $1\leq q\leq \infty$. If 
$$H(\bb \sigma)/d>\left(1/p-1/q\right)_+,$$
then $\mathcal N^{\bb r,\bb \sigma}_{p'}(\L_p([0,1]^d))\subset \L_q([0,1]^d)$ and, for all $s\in\mathcal N^{\bb r,\bb \sigma}_{p'}(\L_p([0,1]^d))$,  $$\|s\|_q \leq  C(d,\bb r,\bb \sigma, p, p',q)  \left( \|s\|_p + N_{\bb r,\bb \sigma,p,p'}(s) \right).$$
\end{prop}
The reader familiar with classical function spaces will have noted the similarity between the definition and the embedding properties of spaces $\mathcal N^{\bb r,\bb \sigma}_{p'}(\L_p([0,1]^d))$ and those of Besov spaces. Before going further, let us recall the definition of the latter according to~\cite{SchmeisserTriebel}, for instance. We denote by $(\textbf{b}_1,\ldots,\textbf{b}_d)$ the canonical basis of $\R^d$ and set $\mathcal{R}=[0,1]^\d$. For all $\bb \sigma=(\sigma_1,\ldots,\sigma_d) \in (0,+\infty)^d$, $0<p,p'\leq \infty$, $s\in\L_p([0,1]^d)$, $h>0$ and $1\leq l\leq d$, we define 
$$\mathcal{R}(\sigma_l,h)=\{x\in [0,1]^d \text{ s.t. } x, x+h\textbf{b}_l,\ldots, x+(\lfloor\sigma_l\rfloor +1)h\textbf{b}_l \in \mathcal{R}\},$$
$$\Delta^{\sigma_l}_{h\textbf{b}_l}s(x)=\sum_{k=0}^{\lfloor\sigma_l\rfloor+1}\binom{\lfloor\sigma_l\rfloor+1}{k}(-1)^{\lfloor\sigma_l\rfloor+1-k}s(x+kh\textbf{b}_l), \text{ for } x\in \mathcal{R}(\sigma_l,h),$$
$$\omega_{\sigma_l}^{(l)}(s,y,\mathcal{R})_p=\sup_{0<h\leq y}\|\Delta^{\sigma_l}_{h\textbf{b}_l}s\BBone_{\mathcal{R}(\sigma_l,h)}\|_p, \text{ for } y\geq 0,$$
\[|s|_{\bb \sigma,p,p'}=
\begin{cases}
\sum_{l=1}^d \left(\int_{0}^\infty \left[y^{-\sigma_l}\omega_{\sigma_l}^{(l)}(s,y,\mathcal{R})_p\right]^{p'}\frac{\d y}{y}\right)^{1/p'}
&\text { if } 0 < p' < \infty\\
\sum_{l=1}^d \left(\sup_{y>0} y^{-\sigma_l}\omega_{\sigma_l}^{(l)}(s,y,\mathcal{R})_p\right)
&\text { if } p' = \infty.
\end{cases}\]
For $\bb \sigma=(\sigma_1,\ldots,\sigma_d) \in (0,+\infty)^d$, $0<p,p'\leq \infty$, we denote by $\mathscr B^{\bb \sigma}_{p'}\left(\L_p([0,1]^d)\right)$ the space of all measurable functions $s\in\L_p([0,1]^d)$ such that $|s|_{\bb \sigma,p,p'}$ is finite. 
According to the proposition below, Besov spaces $\mathscr B^{\bb \sigma}_{p'}\left(\L_p([0,1]^d)\right)$ are embedded in spaces $\mathcal N^{\bb r,\bb \sigma}_{p'}(\L_p([0,1]^d)).$
\begin{prop}\label{Besovsmoothness}
Let $\bb \sigma\in \prod_{l=1}^d (0,r_l+1)$, $0<p<\infty$ and $0<p'\leq \infty$.
For all $s\in \mathscr B^{\bb \sigma}_{p'}\left(\L_p([0,1]^d)\right)$, $$N_{\bb r,\bb \sigma,p,p'}(s) \leq C(d,\bb r,\bb \sigma,p,p') |s|_{\bb \sigma,p,p'}.$$
\end{prop}
We shall not give a proof of that proposition here, since it relies exactly on the same arguments as those used by~\cite{Hochmuth} in the proof of Theorem 4.1 (beginning of page 197) combined with Inequality (14) in the same reference. It should be noticed that the space $\mathcal N^{\bb r,\bb \sigma}_{p'}(\L_p([0,1]^d))$ is in general larger than $\mathscr B^{\bb \sigma}_{p'}(\L_p([0,1]^d))$. Indeed, contrary to $\mathscr B^{\bb \sigma}_{p'}(\L_p([0,1]^d))$, the space $\mathcal N^{\bb r,\bb \sigma}_{p'}(\L_p([0,1]^d))$ contains discontinuous functions (piecewise polynomials, for instance) even for $H(\bb \sigma)/d > 1/p$. 

We are now able to state approximation rates over anisotropic classes of the form
$$\mathcal S(\bb r,\bb\sigma,p,p',R)=\{s\in\L_p([0,1]^d) \text{ s.t. } N_{\bb r,\bb \sigma,p,p'}(s)\leq R \},$$
where $\bb \sigma \in \prod_{l=1}^d (0,r_l+1)$, $0<p,p'\leq \infty$ and $R>0$, thus extending the result of DeVore and Yu~\cite{DeVoreYu} (Corollary 3.3), which is only devoted to functions with isotropic smoothness. The approximation rate is related to the harmonic mean $H(\bb \sigma)$ of $\sigma_1,\ldots,\sigma_d$, which in case of isotropic smoothness of order $\sigma$, \textit{i.e.} if $\sigma_1=\ldots=\sigma_d=\sigma$, reduces to $\sigma$.  
\begin{theo}\label{ApproxAnisotropicNC}
Let $R>0$, $\bb \sigma=(\sigma_1,\ldots,\sigma_d) \in \prod_{l=1}^d (0,r_l+1)$, $0<p < \infty$ and $1\leq q\leq \infty$ such that   
\begin{equation*}
H(\bb \sigma)/d>(1/p-1/q)_+.
\end{equation*}
Assume that $s\in\mathcal S(\bb r,\bb \sigma,p,p',R)$, where $p'=\infty$ if $0<p\leq 1$ or $p\geq q$, and $p'=p$ if $1<p<q$.  
Then, for all $k\in\N$, there exists some partition $m$ of $[0,1]^d$ into dyadic rectangles, that may depend on $s,d,\bb r,\bb \sigma,p$ and $q$, such that  
$$|m|\leq C_1(d,\bb \sigma,p) 2^{kd}$$
and 
\begin{equation}
\inf_{t\in S_{(m,\bb r)}}\|s-t\|_q \leq C_2(d,\bb r,\bb \sigma,p,q) R 2^{-kH(\bb \sigma)}.
\label{eq:boundNC}
\end{equation}
\end{theo}
The same result still holds whatever $0<p'\leq \infty$ if $0<p\leq 1$ or $p\geq q$, and whatever $0<p'\leq p$ if $1<p<q$, as a straightforward consequence of Theorem~\ref{ApproxAnisotropicNC} and Inequality~\eqref{eq:embeddingNN}. Denoting by $\mathcal M_D$, $D\in\N^\star$, the set of all the partitions of $[0,1]^d$ into $D$ dyadic rectangles, we obtain uniform approximation rates simultaneously over a wide range of classes $\mathcal S(\bb r,\bb \sigma,p,p',R)$ by considering the nonlinear approximating space $\cup_{m \in\mathcal{M}_D} S_{(m,\bb r)}$. That property is stated more precisely in Corollary~\ref{uniformNC} below, which can be immediately derived from Theorem~\ref{ApproxAnisotropicNC}. 
\begin{corol}\label{uniformNC}
Let $R>0$, $\bb \sigma=(\sigma_1,\ldots,\sigma_d) \in \prod_{l=1}^d (0,r_l+1)$, $0<p < \infty$, $0<p' \leq\infty$ and $1\leq q\leq \infty$ satisfying the assumptions of Theorem~\ref{ApproxAnisotropicNC}. For all $D\geq C_1(d,\bb \sigma,p)$, where $C_1(d,\bb \sigma,p)$ is given by Theorem~\ref{ApproxAnisotropicNC},
$$\sup_{s\in\mathcal S(\bb r,\bb \sigma,p,p',R)} \inf_{t\in \cup_{m \in\mathcal{M}_D} S_{(m,\bb r)}}\|s-t\|_q \leq C'_2(d,\bb r,\bb \sigma,p,q) R D^{-H(\bb \sigma)/d}.$$
\end{corol} 

\noindent
We also propose of a more refined version of Theorem~\ref{ApproxAnisotropicNC} that allows to take into account constraints on the minimal dimensions of the dyadic rectangles, which will prove most useful for estimation purpose in the next section. We recall that $\underline{\bb \sigma}=\min_{1\leq l\leq d} \sigma_l.$ 
\begin{theo}\label{ApproxAnisotropicC}
Let $J\in\N$, $R>0$, $\bb \sigma=(\sigma_1,\ldots,\sigma_d) \in \prod_{l=1}^d (0,r_l+1)$, $0<p < \infty$,  $0<p'\leq \infty$ and $1\leq q\leq \infty$ such that 
$$H(\bb \sigma)/d>\left(1/p-1/q\right)_+.$$
Assume that $s\in\mathcal S(\bb r,\bb \sigma,p,p',R)$, where $p'=\infty$ if $0<p\leq 1$ or $p\geq q$, and $p'=p$ if $1<p<q$.  
Then, for all $k\in\N$, there exists some partition $m$ of $[0,1]^d$, that may depend on $s,d,\bb r,\bb \sigma,p$ and $q$, only contains dyadic rectangles with sidelength at least $2^{-J\underline{\bb \sigma}/\sigma_l}$ in the $l$-th direction, $l=1,\ldots,d$, and satisfies both  
$$|m|\leq C_1(d,\bb \sigma,p) 2^{kd}$$
and 
\begin{equation}
\inf_{t\in S_{(m,\bb r)}}\|s-t\|_q 
\leq C_3(d,\bb r,\bb \sigma,p,q) R\left(2^{-Jd\left(H(\bb \sigma)/d-(1/p-1/q)_+\right) \underline{\bb \sigma}/H(\bb \sigma)}+2^{-kH(\bb \sigma)}\right).
\label{eq:boundC}
\end{equation}
\end{theo}

\noindent
\textbf{\textit{Remark:}} Given $J\in\N$, that theorem relies on applying the approximation algorithm of Section~\ref{sec:approxalgo} to an approximation of $s$ from $S_{(m_J,\bb r)}$, where $m_J$ is the partition of $[0,1]^d$ into the dyadic rectangles from $\mathcal D^{\bb \sigma}_J.$ Thus, the term $2^{-Jd\left(H(\bb \sigma)/d-(1/p-1/q)_+\right) \underline{\bb \sigma}/H(\bb \sigma)}$ in~\eqref{eq:boundC}, which is of order $(\dim(S_{(m_J,\bb r)}))^{-\left(H(\bb \sigma)/d-(1/p-1/q)_+\right)}$, corresponds with an upper-bound for the linear approximation error $\inf_{t\in S_{(m_J,\bb r)}} ||s-t||_q.$ The upper-bound~\eqref{eq:boundC} is of the same order as~\eqref{eq:boundNC} -- up to a real that only depends on $d,\bb r,\bb \sigma, p,q$ -- as long as
\begin{equation}
k\leq J \frac{\underline{\bb \sigma}}{H(\bb \sigma)}\left(\frac{H(\bb \sigma)}{d}-\left(\frac{1}{p}-\frac{1}{q}\right)_+\right)\frac{d}{H(\bb \sigma)}.
\label{eq:restrict}
\end{equation} 
If $p\geq q$ and $\underline{\bb \sigma}=H(\bb \sigma)$, \textit{i.e.} if $s$ has homogeneous and isotropic smoothness, then that condition simply amounts to $k\leq J$. Otherwise, Condition~\eqref{eq:restrict} is all the more stringent as $p$ is small by comparison with $q$ or as $\underline{\bb \sigma}$ is small by comparison with $H(\bb \sigma)$, \textit{i.e.} all the more stringent as inhomogeneity or anisotropy are pronounced.

Given $J\in\N$, let us denote by $\mathcal{M}_D^J$ the set of all the partitions into $D$ dyadic rectangles with sidelengths $\geq 2^{-J}$, for $D\in\N^\star$. We can still obtain uniform approximation rates simultaneously over a wide range of classes $\mathcal S(\bb r,\bb \sigma,p,p',R)$ under the constraint that the piecewise polynomial approximations are built over dyadic rectangles with sidelengths $\geq 2^{-J}$, by introducing this time the nonlinear approximation space $\cup_{m \in\mathcal{M}^J_D} S_{(m,\bb r)}$. Indeed, as for all $\bb \sigma=(\sigma_1,\ldots,\sigma_d) \in \prod_{l=1}^d (0,r_l+1)$ and $l=1,\ldots,d$, $2^{-J\underline{\bb \sigma}/\sigma_l} \geq 2^{-J}$, a straightforward consequence of Theorem~\ref{ApproxAnisotropicC} is Corollary~\ref{uniformC} below. 
\begin{corol}\label{uniformC}
Let $J\in\N$, $R>0$, $\bb \sigma=(\sigma_1,\ldots,\sigma_d) \in \prod_{l=1}^d (0,r_l+1)$, $0<p < \infty$, $0<p' \leq\infty$ and $1\leq q\leq \infty$ satisfying the assumptions of Theorem~\ref{ApproxAnisotropicC}. For all $D\geq C_1(d,\bb \sigma,p)$, where $C_1(d,\bb \sigma,p)$ is given by Theorem~\ref{ApproxAnisotropicC},
\begin{align*}
&\sup_{s\in\mathcal S(\bb r,\bb \sigma,p,p',R)} \inf_{t\in \cup_{m \in\mathcal{M}^J_D} S_{(m,\bb r)}}\|s-t\|_q\\ 
&\leq C'_3(d,\bb r,\bb \sigma,p,q) R\left(2^{-Jd\left(H(\bb \sigma)/d-(1/p-1/q)_+\right) \underline{\bb \sigma}/H(\bb \sigma)}+2^{-kH(\bb \sigma)}\right).
\end{align*}
\end{corol}

\section{Application to density estimation}\label{sec:estimation}

This section aims at illustrating the interest of the previous approximation results in statistics. More precisely, placing ourselves in the density estimation framework, we show that combining estimation via dyadic piecewise polynomial selection and the aforementioned approximation results leads to a new density estimator which is able to adapt to the unknown smoothness of the function to estimate, even though it is both anisotropic and inhomogeneous. Besides, we explain how such a procedure can be implemented efficiently.

\subsection{Framework and notation}

Let $n\in\N$, $n\geq 4$, we observe independent and identically distributed random variables $Y_1,\ldots,Y_n$ defined on the same measurable space $(\Omega,\mathcal A)$ and taking values in $[0,1]^d$. We assume that $Y_1,\ldots,Y_n$ admit the same density $s$ with respect to the Lebesgue measure $\lambda_d$ on $[0,1]^d$ and that $s\in\L_2([0,1]^d)$. We denote by $P_s$ the joint distribution of $(Y_1,\ldots,Y_n)$, that is the probability measure with density 
$$\frac{dP_s}{d\lambda_d^{\otimes n}} : (y_1,\ldots,y_n) \in [0,1]^d\times \ldots\times[0,1]^d \longmapsto \prod_{i=1}^n s(y_i),$$
while $\P_s$ stands for the underlying probability measure on $(\Omega,\mathcal A)$, so that for all product $B$ of $n$ rectangles of $[0,1]^d$ 
$$P_s(B)=\P_s(\{\omega\in\Omega \text{ s.t. } (Y_1(\omega),\ldots,Y_n(\omega)) \in B\}).$$
The expectation and variance associated with $\P_s$ are denoted by $\E_s$ and $\Var_s$.

\subsection{Dyadic piecewise polynomial estimators}\label{sec:onemodel}

Let $m$ be some partition of $[0,1]^d$ into dyadic rectangles and $\bb \rho=(\bb \rho_K)_{K\in m}$ a sequence such that, for all $K\in m$, $\bb \rho_K=(\rho_K(1),\ldots,\rho_K(d))\in \N^d$. We denote by $S_\mrho$ the space of all functions $t:[0,1]^d\rightarrow \R$ such that, for all $K\in m$, $t$ is polynomial with degree $\leq \rho_K(l)$ in the $l$-th direction on the rectangle $K$. In particular, if $\bb \rho$ is constant and equal to $\bb r$, then $S_\mrho$ coincides with the space $S_\mr$ introduced in Section~\ref{sec:approximation}. 
Let $\langle .,. \rangle$ be the usual scalar product on $\L_2([0,1]^d)$. We recall that $s$ minimizes over $t\in\L_2([0,1]^d)$ 
$$\|s-t\|_2^2-\|s\|_2^2=\|t\|_2^2-2\langle t,s \rangle = \E_s[\gamma(t)],$$
where 
\begin{equation*} 
\gamma(t)=\|t\|_2^2 - \frac{2}{n}\sum_{i=1}^n t(Y_i)
\end{equation*}
only depends on the observed variables. Thus, a natural estimator of $s$ with values in $S_\mrho$ is 
\begin{equation*}
\hat s_\mrho = \argmin{t\in S_\mrho}{\gamma(t)},
\label{eq:onemodel} 
\end{equation*}
that we will call a dyadic piecewise polynomial estimator. Such an estimator is just a projection estimator of $s$ on $S_\mrho$. Indeed, if for each dyadic rectangle $K$ we set $\Lambda(\bb \rho_K)=\prod_{l=1}^d \{0,\ldots,\rho_K(l)\}$ and denote by $(\Phi_{K,\bb k})_{\bb k\in\Lambda(\bb \rho_K)}$ an orthonormal basis of the space of polynomial functions over $K$ with degree $\leq \rho_K(l)$ in the $l$-th direction, then simple computations lead to 
$$\hat s_\mrho = \sum_{K\in m}\sum_{\bb k\in\Lambda(\bb \rho_K)} \left(\frac{1}{n}\sum_{i=1}^n\Phi_{K,\bb k}(Y_i)\right)\Phi_{K,\bb k}.$$
For theoretical reasons, we shall choose in the remaining of the article an orthonormal basis $(\Phi_{K,\bb k})_{\bb k\in\Lambda(\bb \rho_K)}$ derived from the Legendre polynomials in the following way. Let $(Q_j)_{j\in \N}$ be the orthogonal family of the Legendre polynomials in $\L_2([-1,1])$. For  $K=\prod_{l=1}^{d}[u_i,v_i]$ rectangle of $[0,1]^d$, $\bb k=(k(1),\ldots,k(d)) \in \N^d$ and $x=(x_1,\ldots,x_d)\in [0,1]^d$, we set 
$$\pi(\bb k)=\prod_{l=1}^d (2k(l)+1)$$
and
$$\Phi_{K,\bb k}(x) = 
\sqrt{\frac{\pi(\bb k)}{\lambda_{d}(K)}}\prod_{l=1}^{d}Q_{k(l)}\left(\frac{2x_l-u_l-v_l}{v_l-u_l}\right)\BBone_{K}(x).$$
We recall that, for all $j\in\N$,  $Q_j$ satisfies 
$$\|Q_j\|_\infty =1
\quad \text{ and }\quad
\|Q_j\|_2^2=\frac{2}{(2j+1)}.$$
Therefore, for $K$ rectangle in $[0,1]^d$ and $\bb \rho_K \in \N^d$, $(\Phi_{K,\bb k})_{\bb k\in\Lambda(\bb \rho_K)}$ is a basis of the space of piecewise polynomial functions with support $K$ and degree $\leq \rho_K(l)$ in the $l$-th direction, which is orthonormal for the norm $\|.\|_2$ and satisfies 
\begin{equation}
\|\Phi_{K,\bb k}\|^2_\infty = \frac{\pi(\bb k)}{\lambda_d(K)}.
\label{eq:supPhi}
\end{equation}

For each partition $m$ of $[0,1]^d$ into dyadic rectangles and each $\bb \rho=(\bb \rho_K)_{K\in m} \in (\N^d)^{|m|}$, we can evaluate the performance of $\hat s_\mrho$ by giving an upper-bound for its quadratic risk. For that purpose, we introduce the orthogonal projection $s_\mrho$ of $s$ on $S_\mrho$, the dimension $\dim( S_\mrho)$ of $S_\mrho$, \textit{i.e.} 
$$\dim( S_\mrho)=\sum_{K\in m} |\Lambda(\bb \rho_K)|=\sum_{K\in m}\prod_{l=1}^d(\rho_K(l)+1),$$ and define $\bb{\rho_{\max}}=(\rho_{\max}(1),\ldots,\rho_{\max}(d))$ 
by
\begin{equation} 
\rho_{\max}(l)=\max_{K\in m}\rho_K(l), l=1,\ldots,d.
\label{eq:barrho}
\end{equation}
\begin{prop}\label{onemodel}
Let $m$ be a partition of $[0,1]^d$ into dyadic rectangles and $\bb \rho=(\bb \rho_K)_{K\in m} \in (\N^d)^{|m|}$.
If $s\in\L_2([0,1]^d)$, then 
\begin{equation*}
\E_s\left[\|s-\hat s_\mrho\|_2^2\right] = \|s- s_\mrho\|_2^2 + \frac{1}{n} \sum_{K\in m}\sum_{\bb k\in\Lambda(\bb \rho_K)} \Var_s(\Phi_{K,\bb k}(Y_1)).
\label{eq:decomponemodel}
\end{equation*}
If $\|s\|_\infty$ is finite, then
\begin{equation*}
\E_s\left[\|s-\hat s_\mrho\|_2^2\right] \leq \|s- s_\mrho\|_2^2 + \pi(\bb{\rho_{\max}}) \|s\|_\infty \frac{\dim( S_\mrho)}{n}.
\label{eq:upperonemodel}
\end{equation*}
\end{prop} 

\noindent
\textbf{\textit{Proof: }} Pythagoras' Equality gives  
\begin{equation*}
\E_s\left[\|s-\hat s_\mrho\|_2^2\right] = \|s- s_\mrho\|_2^2 + \E_s\left[\|s_\mrho-\hat s_\mrho\|_2^2\right].
\end{equation*}
Then, we deduce the first equality in Proposition~\ref{onemodel} from the expressions of $\hat s_\mrho$ and $s_\mrho$ in the orthonormal basis $(\Phi_{K,\bb k})_{K\in m, \bb k\in\Lambda(\bb \rho_K)}$ of $S_\mrho$ and the fact that $Y_1,\ldots,Y_n$ are independent and identically distributed.

If $s$ is bounded, we deduce from~\eqref{eq:supPhi} that, for all $K\in m$ and $\bb k\in\Lambda(\bb \rho_K)$,
\begin{equation*}
\E_s\left[\Phi^2_{K,\bb k}(Y_1)\right] 
\leq \langle s,\BBone_K \rangle \frac{\pi(\bb k)}{\lambda_d(K)}
\leq \|s\|_\infty \pi(\bb{\rho_{\max}}),
\end{equation*}
hence the upper-bound for $\E_s\left[\|s-\hat s_\mrho\|_2^2\right]$. $\blacksquare$

\smallskip
\noindent
Thus, we recover that, for bounded densities at least, choosing a model $S_\mrho$ that realizes a good compromise between the approximation error and the dimension of the model leads to an estimator $\hat s_\mrho$ with small risk. Such a choice reveals in fact optimal for densities presenting the kind of smoothness described in Section~\ref{sec:approxrate}. More precisely, for $\bb \sigma \in (0,+\infty)^d$, $0<p,p'\leq\infty$, $R>0$ and $L>0$, we set $\lfloor \bb \sigma\rfloor = (\lfloor \sigma_1\rfloor,\ldots,\lfloor \sigma_d\rfloor )$ and consider the class $\mathcal P(\bb \sigma, p, p', R, L)$ of all the probability densities $s$ with respect to $\lambda_d$ such that $s\in\mathcal S(\lfloor \bb \sigma\rfloor + 1, \bb \sigma, p,p',R)$ and $\|s\|_\infty \leq L$. Thanks to the upper-bound of Proposition~\ref{onemodel}, we obtain in Proposition~\ref{minimaxonemodel} below that any statistical procedure which is able to realize approximately $\inf_{m\in\mathcal M, \bb\rho \in \N^d} \E_s\left[\|s-\hat s_\mrho\|_2^2\right]$, where $\mathcal M$ is the collection of all the partitions of $[0,1]^d$ into dyadic rectangles, enjoys adaptivity properties: it also reaches approximately the minimax risk over a wide range of classes $\mathcal P(\bb \sigma, p, p', R, L)$.

\begin{prop}\label{minimaxonemodel}
For $0<p<\infty$, let $p'=\infty$ when $0<p\leq 1$ or $p\geq 2$, and $p'=p$ when $1<p<2$. For all $L>0$ and $R\geq n^{-1/2}$, if $\bb \sigma\in(0,+\infty)^d$ and $0<p<\infty$ satisfy $H(\bb\sigma)/d>(1/p-1/2)_+$, then 
\begin{align*}
&\sup_{s\in \mathcal P(\bb \sigma, p, p', R, L)} \inf_{m\in\mathcal M, \bb\rho \in \N^d} \E_s\left[\|s-\hat s_\mrho\|_2^2\right]\\
&\leq C(d,\bb \sigma, p,L) \left(Rn^{-H(\bb \sigma)/d}\right)^{2d/(d+2H(\bb \sigma))}\\
&\leq  C(d,\bb \sigma, p,L) \inf_{\hat s}\sup_{s\in \mathcal P(\bb \sigma, p, p', R, L)} \E_s\left[\|s-\hat s\|_2^2\right]
\end{align*}
where the last infimum is taken over all the estimators $\hat s$ of $s$.
\end{prop}
\noindent
\textbf{\textit{Proof: }} Let us fix $\bb \sigma, p, p', R, L$ satisfying the assumptions of Proposition~\ref{minimaxonemodel} and choose $\bb \rho=\lfloor \bb \sigma \rfloor +1$. For all $s\in\mathcal P(\bb \sigma, p, p', R, L)$, we deduce from Proposition~\ref{onemodel} and Theorem~\ref{ApproxAnisotropicNC} that
\begin{equation*}
\inf_{m\in\mathcal M, \bb\rho \in \N^d} \E_s\left[\|s-\hat s_\mrho\|_2^2\right]
\leq C(d, \bb \sigma,p,L)\inf_{k\in\N}  \left\{R^2 2^{-2kH(\bb \sigma)} +  \frac{2^{kd}}{n}\right\}.
\end{equation*}
We then choose $k_\star$ as the greatest integer $k\in\N$ such that $2^{kd}/n \leq R^2 2^{-2kH(\bb \sigma)}$, \textit{i.e.} such that $2^{k}\leq (nR^2)^{1/(d+2H(\bb \sigma))}$ so as to bound the infimum on the right-hand side, which provides the first inequality in Proposition~\ref{minimaxonemodel}.

Let us define the Besov class $\mathcal B(\bb \sigma, p, p', R, L)$ of all the probability densities $s$ with respect to $\lambda_d$ such that 
$|s|_{\bb \sigma, p,p'} \leq R$ (where $|.|_{\bb \sigma, p,p'}$ is defined in Section~\ref{sec:approxrate}) and $\|s\|_\infty \leq L$. We deduce from Proposition~\ref{Besovsmoothness} that, for $0<p\leq 1$ or $p\geq 2$, there exists some positive real $C(\bb \sigma, p)$ such that $\mathcal P(\bb \sigma, p, \infty, R, L)$ contains $\mathcal B(\bb \sigma, \infty, \infty, C(\bb \sigma, p)R, L)$, and, for $1<p<2$,  there exists some positive real $C(\bb \sigma, p)$ such that $\mathcal P(\bb \sigma, p, p, R, L)$ contains $\mathcal B(\bb \sigma, p, p, C(\bb \sigma, p)R, L)$. Besides, according to Triebel~\cite{TriebelArticle} (Proposition 10), for all $\epsilon>0$, the Kolmogorov $\epsilon$-entropy in $\L_2([0,1]^d)$ of the Besov space $\mathscr B^{\bb \sigma}_q(\L_q([0,1]^d))$ is $\epsilon^{-H(\bb\sigma)/d}$ for $H(\bb\sigma)/d>(1/q-1/2)_+$. Thus, the second inequality in Proposition~\ref{minimaxonemodel} follows from the lower-bounds for minimax risks proved in~\cite{YangBarron} (Proposition 1, $ii)$). $\blacksquare$

In the sequel, our problem will thus be to build a statistical procedure that requires no prior knowldege on $s$ but whose risk behaves almost as $\inf_{m\in\mathcal M, \bb\rho \in \N^d} \E_s\left[\|s-\hat s_\mrho\|_2^2\right]$.

\subsection{Dyadic piecewise polynomial selection}
Let us fix $\bb r_\star\in\N^d$, $J_\star\in\N$, and denote by $\M$ the set of all partitions of $[0,1]^d$ into dyadic rectangles with sidelengths at least $2^{-J_\star}$. We consider the family $\Mdeg$ of all couples $(m,\bb \rho)$ with $m\in\mathcal{M}_\star$ and $\bb \rho=(\bb \rho_K)_{K\in m} $ such that, for all $K\in m$, $\bb\rho_K\in \Lambda(\bb{r_\star})$. Ideally, we would like to choose the couple $(m,\bb \rho)$ that minimizes $\E_s\left[\|s-\hat s_{(m,\bb \rho)}\|_2^2\right]$ among the elements of $\Mdeg$. This is hopeless without knowing $s$, but from Pythagora's Equality and Proposition~\eqref{eq:onemodel}, we have, for all $\mrho\in\Mdeg$, 
$$\E_s\left[\|s-\hat s_\mrho\|_2^2\right] - \|s\|_2^2= - \|s_\mrho\|_2^2 + \frac{1}{n} \sum_{K\in m}\sum_{\bb k\in\Lambda(\bb \rho_K)} \Var_s(\Phi_{K,\bb k}(Y_1)).$$ 
Thus, we propose to select an adequate partition $\hat m$ and the associated sequence of maximal degrees $\bb{\hat\rho} =(\bb{\hat\rho}_K)_{K\in\hat m}$ from the data so that 
\begin{align*}
(\hat m,\bb{\hat\rho}) 
&= \argmin{\mrho \in \Mdeg}{\{- \|\hat s_\mrho\|_2^2 + \pen\mrho\}}\\
&=\argmin{\mrho \in \Mdeg}{\{\gamma(\hat s_\mrho) + \pen\mrho\}}
\end{align*}  
where $\pen:\Mdeg\rightarrow\R^+$ is a so-called penalty function. We then estimate the density $s$ by 
\begin{equation*}
\tilde s = \hat s_\hatmrho.
\label{eq:stilde}
\end{equation*}  
According to the proof of Proposition~\ref{minimaxonemodel}, in view of proving the adaptivity of the penalized estimator $\tilde s$, the penalty $\pen$ should be chosen so that $\tilde s$ satisfies an inequality akin to
\begin{equation}
\E_s[\|s-\tilde s\|^2_2] \leq C \min_{\mrho \in \Mdeg} \left\{\|s- s_\mrho\|_2^2 + \frac{\dim(S_\mrho)}{n}\right\}
\label{eq:goal}
\end{equation}
where $C$ is a positive real that does not depend on $n$.

In order to define an adequate form of penalty, we introduce the set $\mathcal{D}_\star$ of all dyadic rectangles of $[0,1]^d$ with sidelengths $\geq 2^{-J_\star}$ and, for all $K\in\mathcal D_\star$ and $\bb k\in\Lambda(\bb r_\star)$, we set 
$$\hat \sigma^2_{K,\bb k}=\frac{1}{n(n-1)}\sum_{i=2}^n\sum_{j=1}^{i-1} \left(\Phi_{K,\bb k}(Y_i)-\Phi_{K,\bb k}(Y_j)\right)^2,$$
which is an unbiased estimator of $\Var_s(\Phi_{K,\bb k}(Y_1))$. We also set
\begin{equation*}  
\widehat M_{1,\star}=\frac{1}{n}\max_{K\in\mathcal D_\star} \sum_{\bb k\in\Lambda(\bb r_\star)} \sqrt{\frac{\pi(\bb k)}{\lambda_d(K)}}\left|\sum_{i=1}^n \Phi_{K,\bb k}(Y_i)\right|
\quad\text{and}\quad
\widehat M_{2,\star} = \frac{1}{n} \max_{K\in\mathcal D_\star} \max_{\bb k\in\Lambda(\bb {r_\star})} \sum_{i=1}^n \Phi^2_{K,\bb k} (Y_i),
\end{equation*}
that overestimate respectively 
\begin{equation*}  
\max_{\mrho \in \Mdeg} \|s_\mrho\|_\infty
\quad\text{and}\quad
\max_{K\in \mathcal D_\star} \max_{\bb k\in\Lambda\left(\bb{r_\star}\right)} \E_s\left[\Phi^2_{K,\bb k} (Y_1)\right].
\end{equation*}
The following theorem suggests a form of penalty yielding an inequality close to~\eqref{eq:goal}.
\begin{theo}\label{selection}
Let $\bb {r_\star}\in\N^d$ and $J_\star\in\N$ be such that $|\Lambda(\bb {r_\star})|\leq \max\{\exp(n)/n, n^d\}$ and $2^{dJ_\star} \leq n/\log(n|\Lambda(\bb {r_\star})|)$. Let $(L_{\mrho})_{\mrho\in\Mdeg}$ be a family of nonnegative real numbers, that may depend on $n$, satisfying 
\begin{equation}
\sum_{\mrho\in\Mdeg} \exp(-L_{\mrho}|m|) \leq 1.
\label{eq:Sigma}
\end{equation}
If $s$ is bounded and $\pen$ is defined on $\Mdeg$ by 
\begin{align*}
\pen\mrho
&=\frac{1}{n}\sum_{K\in m} \sum_{\bb k\in\Lambda(\bb{\rho}_K)}\left(\kappa_1 \hat\sigma^2_{K,\bb k} + \kappa_2 \pi(\bb k)\right)\\
&+\left(\left(\kappa_3\widehat M_{2,\star} +\kappa_4\pi(\bb r_\star)\right)|\Lambda(\bb r_\star)|+\kappa_5\widehat M_{1,\star}\right)\frac{L_{\mrho}|m|}{n} 
\end{align*}
where $\kappa_1,\ldots,\kappa_5$ are large enough positive constants, then
\begin{align*}
\E_s\left[\|s-\tilde s\|_2^2\right]
&\leq \min_{\mrho\in \Mdeg}\Bigg\{\kappa'_1 \|s-s_\mrho\|_2^2 
+\kappa'_2\frac{1}{n} \sum_{K\in m}\sum_{\bb k\in\Lambda(\bb \rho_K)} \Var_s(\Phi_{K,\bb k}(Y_1))\\
&+\kappa'_3\pi(\bb r_\star)\frac{\dim(S_\mrho)}{n}
+\kappa'_4\pi(\bb r_\star)|\Lambda(\bb r_\star)|\|s\|_\infty\frac{L_{\mrho}|m|}{n}\Bigg\}\notag\\
&+\kappa'_5 \|s\|^2_\infty\pi(\bb r_\star)|\Lambda(\bb r_\star)| \frac{1}{n}.
\end{align*}
where $\kappa'_1,\ldots,\kappa'_5$ are positive reals, $\kappa'_1,\ldots,\kappa'_4$ only depend on $\kappa_1,\ldots,\kappa_5$, and $\kappa'_5$ also depends on $d$.
\end{theo}
Thus, the penalty associated to each $\mrho\in\Mdeg$ is composed of two terms: an additive term that overestimates the variance over the model $S_\mrho$, and a term linear in the size of the partition $m$, up to the weight $L_\mrho$, that overestimates the upper-bound given in Proposition~\ref{onemodel} for the variance over $S_\mrho$. There remains to choose those weights under the constraint~\eqref{eq:Sigma}. According to Proposition~\ref{choiceL} below, each model in $\Mdeg$ can be assigned the same weight that only depends on $d$ and $\bb r_\star$.

\begin{prop}\label{choiceL}
If $\kappa_1,\ldots,\kappa_5$ are large enough positive constants, then the penalty defined on $\Mdeg$ by 
\begin{align}
\pen\mrho
&=\frac{1}{n}\sum_{K\in m} \sum_{\bb k\in\Lambda(\bb{\rho}_K)}\left(\kappa_1 \hat\sigma^2_{K,\bb k} + \kappa_2 \pi(\bb k)\right)\notag\\
&+\left(\left(\kappa_3\widehat M_{2,\star} +\kappa_4\pi(\bb r_\star)\right)|\Lambda(\bb r_\star)|+\kappa_5\widehat M_{1,\star}\right)\frac{\log(8d|\Lambda(\bb{r_\star})|)|m|}{n} 
\label{eq:pendyadic}
\end{align}
satisfies the assumptions of Theorem~\ref{selection}. Moreover, if $|\Lambda(\bb {r_\star})|\leq \max\{\exp(n)/n, n^d\}$, $2^{dJ_\star} \leq n/\log(n|\Lambda(\bb {r_\star})|)$ and $s$ is bounded, then for such a penalty 
\begin{equation*}
\E_s\left[\|s-\tilde s\|_2^2\right]
\leq \kappa'' \min_{\mrho\in \Mdeg}\left\{\|s-s_\mrho\|_2^2+\pi(\bb r_\star)|\Lambda(\bb r_\star)|\|s\|^2_\infty\frac{\log(8ed|\Lambda(\bb{r_\star})|)|m|}{n}\right\}.
\end{equation*}
where $\kappa''$ is a positive real that only depends on $\kappa_1,\ldots,\kappa_5$ and $d$.
\end{prop}

\begin{figure}[h]
\setlength{\unitlength}{0.4cm}
\begin{center}
\begin{picture}(8,8)
\put(0,0){\line(1,0){8}}
\put(0,0){\line(0,1){8}}
\put(0,8){\line(1,0){8}}
\put(8,8){\line(0,-1){8}}
{{\put(0,4){\line(1,0){8}}}}
{{{\put(4,0){\line(0,1){4}}}}}
{{{\put(0,6){\line(1,0){8}}}}}
{{{\put(4,2){\line(1,0){4}}}}}
\end{picture}
\end{center}
\begin{center}
\pstree[treesep=0.4cm,treefit=tight,levelsep=1.2cm]
 {\TR{$[0,1]\times [0,1]$}}{        
\pstree{\TR{$[0,1]\times\left[0,\frac{1}{2}\right]$}^{(2)}} {           
\TR{$\left[0,\frac{1}{2}\right]\times\left[0,\frac{1}{2}\right]$}^{(1)}
\pstree{\TR{$\left(\frac{1}{2},1\right]\times\left[0,\frac{1}{2}\right]$}}
{\TR{$\left(\frac{1}{2},1\right]\times \left[0,\frac{1}{4}\right]$}^{(2)} 
\TR{$\left(\frac{1}{2},1\right]\times\left(\frac{1}{4},\frac{1}{2}\right]$}}
}
 \pstree{\TR{$[0,1]\times\left(\frac{1}{2},1\right]$}}
 {\TR{$[0,1]\times\left(\frac{1}{2},\frac{3}{4}\right]$}^{(2)} 
 \TR{$[0,1]\times\left(\frac{3}{4},1\right]$}}                
}
\end{center}
\caption{Top: Partition of $[0,1]^2$ into dyadic rectangles. Bottom: Binary tree labeled with the sequence of cutting directions $(2,1,2,2)$ corresponding with that partition.}
\label{flottants:treepart:figure}        
\end{figure}
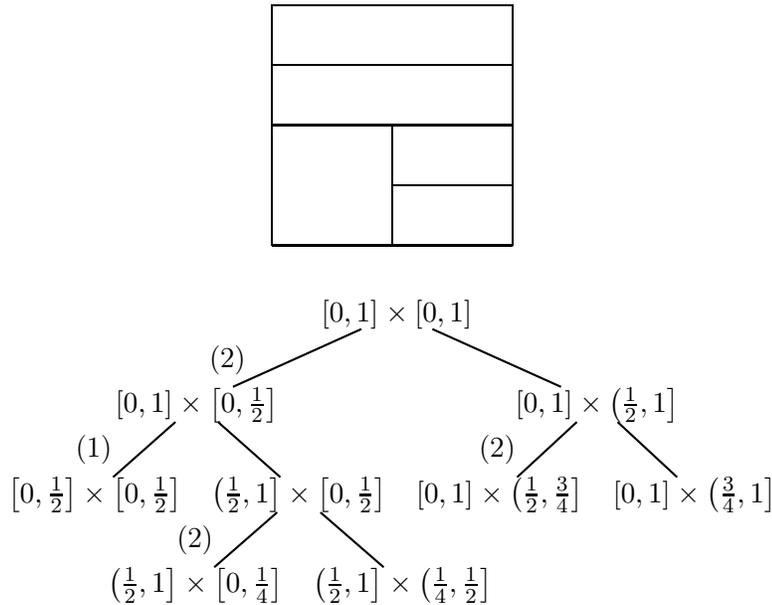

\noindent
\textbf{\textit{Proof:}} First, for all $D\in\N^\star$, the number of partitions of $[0,1]^d$ into $D$ dyadic rectangles satisfies 
\begin{equation}
|\mathcal M_D| \leq (4d)^D.
\label{eq:combinatorics}
\end{equation}
Indeed, as illustrated by Figure~\ref{flottants:treepart:figure}, each partition in $\mathcal M_D$ can be described by a complete dyadic tree with $D$ leaves whose edges are labeled with a sequence of $D-1$ integers in $\{1,\ldots,d\}$ giving the cutting directions to obtain the partition from the unit square. 

\noindent
The number of complete dyadic trees with $D$ leaves is given by the Catalan number  $$\frac{1}{D}\binom{2(D-1)}{D-1}\leq 4^D,$$
hence~\eqref{eq:combinatorics}. We deduce from~\eqref{eq:combinatorics} that, for all positive real $L$,
\begin{align*}
\sum_{\mrho\in\Mdeg} \exp(-L|m|) 
&\leq \sum_{D\in\N^\star} \sum_{m\in \mathcal M_D} \sum_{\bb{\rho} \in \Lambda(\bb{r_\star})^D} \exp(-L|m|)\\
&\leq \sum_{D\in\N^\star} (4d|\Lambda(\bb{r_\star})|)^D \exp(-LD)\\
&\leq 1/\left(\exp\left(L-\log(4d|\Lambda(\bb{r_\star})|)\right)-1\right)
\end{align*}
So, we can choose $L \geq \log(8d|\Lambda(\bb{r_\star})|)$ for Condition~\eqref{eq:Sigma} to be fulfilled. 

Since $\|s\|_\infty \geq 1$, the upper-bound for $\E_s\left[\|s-\tilde s\|_2^2\right]$ is then a straightforward consequence of Theorem~\ref{selection}. $\blacksquare$

\noindent
It is worth pointing out that penalty~\eqref{eq:pendyadic} is more refined than the penalties proposed by~\cite{Klemela} or~\cite{AkakpoDurot} for density estimation via dyadic histogram selection based on a least-squares type criterion. Indeed, when $\bb r_\star$ is null, penalty~\eqref{eq:pendyadic} is not simply proportional to the dimension of the partition.

With a penalty chosen as above, we recover an inequality close to~\eqref{eq:goal}, that allows to prove the adaptivity of $\tilde s$ over a wide range of classes $\mathcal P(\bb \sigma,p,p',R,L)$ as defined in Section~\ref{sec:onemodel}. For that purpose, we introduce 
$$q(d,\bb \sigma,p)=\frac{\underline{\bb \sigma}}{H(\bb \sigma)}\frac{d+2H(\bb \sigma)}{H(\bb \sigma)}\left(\frac{H(\bb \sigma)}{d}-\left(\frac{1}{p}-\frac{1}{2}\right)_+\right)$$
and 
$$w(\bb r_\star) = \pi(\bb r_\star) |\Lambda(\bb r_\star)| \log\left(8ed|\Lambda(\bb r_\star)|\right).$$

\begin{theo}\label{adaptselect}
Let $\bb r_\star\in\N^d$ and $J_\star\in\N$ be such that $|\Lambda(\bb {r_\star})|\leq \max\{\exp(n)/n, n^d\}$ and $J_\star=\max\{ J\in \N \text{ s.t. } 2^{Jd} \leq n/\log(n|\Lambda(\bb r_\star)|)\}$, and $\pen$ be the penalty given by Proposition~\ref{choiceL}. For all $p>0$, let $p'=\infty$ if $0<p\leq 1$ or $p\geq 2$, and $p'=p$ if $1<p<2$. For all $L>0$, $\bb \sigma \in \prod_{l=1}^d (0, r_\star(l)+1)$, $p>0$ such that $H(\bb \sigma)/d >(1/p-1/2)_+$ and $q(d,\bb \sigma, p) >1$, for all $R$ such that $ w(\bb r_\star)/n\leq R^2 \leq (n/\log(n|\Lambda(\bb r_\star)|))^{q(d,\bb \sigma,p)-1},$
\begin{equation*}
\sup_{s\in\mathcal P(\bb \sigma, p,p',R,L)} \E_s\left[\|s-\tilde s\|_2^2\right]
\leq C w(\bb r_\star)^{2H(\bb \sigma)/(d+2H(\bb \sigma))}  \inf_{\hat s}\sup_{s\in \mathcal P(\bb \sigma, p, p', R, L)} \E_s\left[\|s-\hat s\|_2^2\right],
\end{equation*}
where $C$ only depends on $d,\bb \sigma, p,L$ and the penalty constants $\kappa_1,\ldots, \kappa_5$ and the above infimum is taken over all the estimators of $s$.
\end{theo}

\noindent
Thus, if $\bb r_\star$ is chosen as a constant with respect to $n$, then $\tilde s$ reaches the minimax risk, up to a constant factor, over a wide range of classes that contain functions with possibly anisotropic and inhomogeneous smoothness limited by the maximal degrees $\bb r_\star$. Another strategy consists in allowing the maximal degrees $\bb r_\star$ to increase with the sample size $n$, while $w(\bb r_\star)$ varies slowly with $n$. For instance, with $r_\star(l)=\log(n)$ for all $l=1,\ldots,d$, our estimator $\tilde s$ still approximately reaches the minimax risk over a range of classes all the wider as $n$ increases. The price to pay is only a logarithmic factor, proportional to $(\log(\log(n))\log^{2d}(n))^{2H(\bb \sigma)/(d+2H(\bb \sigma))}$ over classes with smoothness $H(\bb \sigma)$. Thus, such a result may be seen as a nonasymptotic and multivariate counterpart of Theorem 1 in Willett and Nowak~\cite{WillettNowak}.

\smallskip
\textbf{\textit{Remark:}} Contrary to~\cite{NeumannVonSachs,Neumann,KLP,Klemela}, we have chosen here the smoothing parameter $J_\star$ independently of the smoothness of $s$, hence the restriction on $q(d,\bb \sigma,p)$, that could disappear otherwise. Setting $\mu_{\bb \sigma}=H(\bb\sigma)/\underline{\bb\sigma}$, the condition $q(d,\bb \sigma,p)>1$ is equivalent to $H(\bb \sigma)/d > \nu(\bb \sigma, p)$, where 
$$\nu(\bb \sigma, p) = \frac{1}{2} \left(\frac{1}{2}\left(\mu_{\bb \sigma}-1\right) + \left(\frac{1}{p}-\frac{1}{2}\right)_+ + \sqrt{\left(\frac{1}{2}\left(\mu_{\bb \sigma}-1\right) + \left(\frac{1}{p}-\frac{1}{2}\right)_+\right)^2+2\left(\frac{1}{p}-\frac{1}{2}\right)_+}\right).$$
In case of isotropic and homogeneous smoothness, \textit{i.e.} when $\mu_{\bb \sigma}=1$ and $p\geq 2$, $q(d,\bb \sigma,p)>1$ is simply equivalent to $H(\bb \sigma)/d>0$. In case of isotropic and inhomogeneous smoothness, \textit{i.e.} when $\mu_{\bb \sigma}=1$ and $p<2$, $q(d,\bb \sigma,p)>1$ is equivalent to  $H(\bb \sigma)/d > \nu(\bb \sigma, p)$ where $ \nu(\bb \sigma, p)\in(1/p-1/2,1/p)$. This is slightly stronger than $H(\bb \sigma)/d >1/p -1/2$, but still better than the restriction $H(\bb \sigma)/d >1/p$ which is often encountered in the literature. Otherwise, $\nu(\bb \sigma,p)$ increases with $\mu_{\bb \sigma}$ and $1/p$, \textit{i.e.} with the anisotropy and the inhomogeneity.

\subsection{Implementing the dyadic piecewise polynomial selection procedure}

We end this article with a brief discussion about the implementation of our estimator $\tilde s$ for the penalty defined in Proposition~\ref{choiceL}. Let us fix the penalty constants $\kappa_1,\ldots,\kappa_5$ and set, for all dyadic rectangle $K\in\mathcal D_\star$ and all $\bb r\in \Lambda(\bb r_\star)$, 
\begin{align*}
\widehat W(K,\bb r)
&=\sum_{\bb k\in\Lambda(\bb r)} \left(-\left(\frac{1}{n}\sum_{i=1}^n\Phi_{K,\bb k}(Y_i)\right)^2 
+\kappa_1 \frac{\hat\sigma^2_{K,\bb k}}{n} + \kappa_2 \frac{\pi(\bb k)}{n}\right)\\ 
&+\frac{\log(8d|\Lambda(\bb{r_\star})|)}{n}\left(\left(\kappa_3\widehat M_{2,\star} +\kappa_4\pi(\bb r_\star)\right)|\Lambda(\bb r_\star)|+\kappa_5\widehat M_{1,\star}\right) 
\end{align*}
and 
$$\hat{\bb r}_K = \argmin{\bb r \in \Lambda(\bb r_\star)} \widehat W(K,\bb r).$$
Given the decomposition of $\hat s_\mrho$ in the basis $(\Phi_{K,\bb k})_{K\in m, \bb k\in \Lambda(\bb\rho_K)}$, the model $\hatmrho$ to select in $\Mdeg$ is characterized by 
$$ \hatmrho = \argmin{\mrho\in\Mdeg} \sum_{K\in m}\widehat W(K,\bb \rho_K),$$ 
so 
$$\hat m=\argmin{m\in\M} \sum_{K\in m} \widehat W(K,\hat{\bb r}_K) \text{ and,  for all } K\in \hat m, \hat{\bb\rho}_K=\hat{\bb r}_K.$$
Thus, the steps leading to $\tilde s$ are 
\begin{enumerate}
\item Compute $\widehat M_{1,\star}$ and $\widehat M_{2,\star}$.
\item For all $K\in\mathcal D_\star$ and all $\bb k \in \Lambda(\bb r_\star)$, compute $\hat \sigma^2_{K,\bb k}$.
\item For all $K\in\mathcal D_\star$, compute $\hat{\bb r}_K$ and $\widehat W(K,\hat{\bb r}_K)$.
\item Determine the best partition $\hat m =\argmintext{m\in\M} \sum_{K\in m} \widehat W(K,\hat{\bb r}_K).$
\item Set, for all $K\in \hat m, \hat{\bb\rho}_K=\hat{\bb r}_K.$ 
\item Compute $\tilde s=\hat s_{\hatmrho}.$
\end{enumerate}
Since $\hat m$ is the partition in $\M$ that minimizes a given additive criterion, it can be determined via the algorithm inspired from Donoho~\cite{DonohoCart} and described in~\cite{BSR} (beginning of Section 3), with a computational complexity at most of order $\mathcal O(|\mathcal D_\star|)$. Therefore, one easily verifies that the whole steps only require a computational complexity at most of order  $\mathcal O(|\Lambda(\bb r_\star)||\mathcal D_\star|)$. Since $|\mathcal D_\star|=(2^{J_\star+1}-1)^d$, if we choose $J_\star$ as prescribed by Theorem~\ref{selection}, then determining $\tilde s$ requires at most $\mathcal{O}(n)$ computations when $\bb r_\star$ is constant, and at most $\mathcal{O}(n\log^d(n))$ when $r_\star(l)=\log(n)$ for all $l=1,\ldots,d$. Last, regarding the choice of the penalty constants $\kappa_1,\ldots,\kappa_5$, they can be calibrated via simulations over a wide collection of test densities. Such a method has already proved to yield good results in practice, even though several constants have to be chosen, as shown for instance in~\cite{ComteRozenholc}.


\section{Proofs of the approximation results}\label{sec:approxproof}

For $j\in\N$ and $K\in\mathcal{D}^{\bb \sigma}_j$, we recall that the children of $K$ are all the dyadic rectangles of $\mathcal{D}^{\bb \sigma}_{j+1}$ that are included in $K$. We will also refer to $K$ as the parent of its children and will often use the fact that the children of $K$ form a partition of $K$ into 
\begin{equation}
\prod_{l=1}^d 2^{\lfloor (j+1)\underline{\bb \sigma}/\sigma_l \rfloor-\lfloor j\underline{\bb \sigma}/\sigma_l \rfloor}
\leq 2^d 2^{d\underline{\bb \sigma}/H(\bb \sigma)}
\label{eq:children}
\end{equation}
dyadic rectangles from $\mathcal{D}^{\bb \sigma}_{j+1}$. 

In all the proofs, the notation $C(\theta)$ stands for a positive real that only depends on the parameter $\theta$, and whose value is allowed to change from one occurrence to another.

\subsection{Proof of Proposition~\ref{embeddingLq}}

For $p\geq q$, Proposition~\ref{embeddingLq} follows from the continuous embedding of $\L_p([0,1]^d)$ in $\L_q([0,1]^d)$. For $p<q$, it corresponds with the second point in the more general result below. 

\begin{prop}\label{changenormaniso}
Let $R>0$, $\bb \sigma=(\sigma_1,\ldots,\sigma_d) \in \prod_{l=1}^d (0,r_l+1)$, $1\leq q\leq \infty$ and $0<p<q$ such that $H(\bb \sigma)/d>1/p-1/q$. For $s\in\L_p([0,1]^d)$, $k\in\N$, and any dyadic rectangle $K\in\mathcal{D}^{\bb\sigma}_k$, we set
\begin{equation*}
e_{\bb r,\bb \sigma, p,k}(s,K)= \inf_{P\in\Pi^{\bb r, \bb \sigma}_{k}} \|(s-P)\BBone_K\|_p.
\end{equation*}
If $N_{\bb r,\bb \sigma, p,\infty}(s)\leq R$, then
\begin{enumerate}[$i)$] 
\item for all $j\in\N$ and all $K\in\mathcal D^{\bb \sigma}_j$,
\begin{equation}
\mathcal{E}_{\bb r}(s,K)_q
\leq C(d,\bb r,\bb \sigma,p,q) \sum_{k\geq j} 2^{-kd(H(\bb \sigma)/d+1/q-1/p)\underline{\bb \sigma} /H(\sigma)} 2^{k\underline{\bb \sigma}} e_{\bb r,\bb \sigma, p,k}(s,K).
\label{eq:changenormrect}
\end{equation}
\item $s\in\L_q([0,1]^d)$ and $\|s\|_q \leq C(d,\bb r,\bb \sigma,p,q) (\|s\|_p+R).$
\end{enumerate} 
\end{prop}

\noindent
\textbf{\emph{Proof: }} Let us fix $1\leq q\leq \infty$, $0<p <q,$  $j\in\N$ and $K\in\mathcal{D}_j^{\bb \sigma}$. For all $k\geq j$, we denote by $\mathcal{C}_k(K)$ the set of all rectangles from $\mathcal{D}_k^{\bb \sigma}$ that are included in $K$. Thus, $\mathcal{C}_j(K)$ is reduced to $\{K\}$,  $\mathcal{C}_{j+1}(K)$ is the set of all the children of $K$, \textit{etc}$\ldots.$ For any rectangle $I\subset [0,1]^d$, we denote by $P_I(s)$ a polynomial function on $I$ with degree $\leq r_l$ in the $l$-th direction such that
$$\|(s-P_I(s))\BBone_I\|_p = \mathcal{E}_{\bb r}(s,I)_p,$$
where $ \mathcal{E}_{\bb r}(s,I)_p$ is defined as in~\eqref{eq:criterion}.
For all $k\geq j$, we set 
$$\Sigma_k(s,K)=\sum_{I\in\mathcal{C}_k(K)} P_I(s) \BBone_I$$
and, in order to alleviate the notation, we simply write $e_k(s,K)$ instead of $e_{\bb r,\bb \sigma, p,k}(s,K)$ in the whole proof.
It should be noticed that $e_k(s,[0,1]^d)=e_{\bb r,\bb \sigma,p,k}(s)$ as defined by~\eqref{eq:pperrork}, and that 
\begin{equation}
e_k(s,K)=\|(s-\Sigma_k(s,K))\BBone_K\|_p =\left(\sum_{I\in\mathcal{C}_k(K)} \mathcal{E}_{\bb r}^p(s,I)_p\right)^{1/p}.
\label{eq:ppapprox}
\notag
\end{equation}
Therefore, 
\begin{equation*}
\|(s-\Sigma_k(s,K))\BBone_K\|_p
\leq \left(\sum_{I\in\mathcal{D}^{\bb \sigma}_k} \mathcal{E}_{\bb r}^p(s,I)_p\right)^{1/p}
= e_{\bb r,\bb \sigma,p,k}(s) \leq 2^{-k\underline{\bb \sigma}}R
\end{equation*}
so that the sequence $(\Sigma_k(s,K))_{k\geq j}$ converges to $s\BBone_K$ in $\L_p([0,1]^d).$ 

Let us prove that $(\Sigma_k(s,K))_{k\geq j}$ also converges to $s\BBone_K$ in $\L_q([0,1]^d).$ We now fix $k\geq j$. When $0<p<q\leq \infty$ as assumed here, Markov Inequality for polynomials asserts that, for all rectangle $I$ of $[0,1]^d$, and all polynomial function $P\in\mathscr{P}_{\bb r}$,
\begin{equation}
\|P\BBone_I\|_q \leq C(d,\bb r,p,q) (\lambda_d(I))^{(1/q-1/p)} \|P\BBone_I\|_p.
\label{eq:Markov}
\end{equation}
We refer to Lemma 5.1 in~\cite{Hochmuth} for a proof (that still holds for $q=\infty$). Let us first assume that $0<p<q<\infty$. We then deduce from~\eqref{eq:Markov} that
\begin{align}
&\|\Sigma_{k+1}(s,K)-\Sigma_k(s,K)\|_q^q\notag\\
&= \sum_{I\in\mathcal{C}_{k+1}(K)}\|(\Sigma_{k+1}(s,K)-\Sigma_k(s,K))\BBone_I\|_q^q\notag\\
&\leq C(d,\bb r,p,q) \sum_{I\in\mathcal{C}_{k+1}(K)} (\lambda_d(I))^{q(1/q-1/p)} \|(\Sigma_{k+1}(s,K)-\Sigma_k(s,K))\BBone_I\|_p^q\notag\\
&\leq C(d,\bb r,p,q) 2^{q(k+1)d(1/p-1/q)\underline{\bb \sigma} /H(\sigma)} \sum_{I\in\mathcal{C}_{k+1}(K)} \|(\Sigma_{k+1}(s,K)-\Sigma_k(s,K))\BBone_I\|_p^q.
\label{eq:decompLq}
\end{align}
Let us also fix $I\in \mathcal{C}_{k+1}(K)$. Then
$$(\Sigma_{k+1}(s,K)-\Sigma_k(s,K))\BBone_I = (P_I(s)-P_{\tilde I}(s))\BBone_I$$
where $\tilde I \in\mathcal{C}_k(K)$ is the parent of $I$. Let $\kappa(p)=2^{1/p}$ if $p<1$, and $\kappa(p)=1$ otherwise. From the (quasi-)triangle inequality satisfied by $\|.\|_p$, we then get 
\begin{align*}
\|(\Sigma_{k+1}(s,K)-\Sigma_k(s,K))\BBone_I\|_p
&\leq \kappa(p) \left(\|(s-P_I(s))\BBone_I\|_p+\|(s-P_{\tilde I}(s))\BBone_I\|_p\right)\\
&\leq \kappa(p) \left(\mathcal{E}_{\bb r}(s,I)_p + \mathcal{E}_{\bb r}(s,\tilde I)_p\right),
\end{align*}
hence, by convexity of $x\mapsto x^q,$
$$\|(\Sigma_{k+1}(s,K)-\Sigma_k(s,K))\BBone_I\|^q_p
\leq 2^{q-1}\kappa^q(p) \left(\mathcal{E}^q_{\bb r}(s,I)_p + \mathcal{E}^q_{\bb r}(s,\tilde I)_p\right).$$
By grouping all the rectangles $I\in\mathcal{C}_{k+1}(K)$ that have the same parent, we obtain
\begin{eqnarray*}
&\sum_{I\in\mathcal{C}_{k+1}(K)}\|(\Sigma_{k+1}(s,K)-\Sigma_k(s,K))\BBone_I\|^q_p \\
&\leq 2^{q-1}\kappa^q(p) \left(\sum_{I\in\mathcal{C}_{k+1}(K)}\mathcal{E}^q_{\bb r}(s,I)_p + 2^{d(1+\underline{\bb \sigma}/H(\bb \sigma))}\sum_{\tilde I\in\mathcal{C}_{k}(K)}\mathcal{E}^q_{\bb r}(s,\tilde I)_p\right).
\end{eqnarray*}
The classical inequality between $\ell_p$ and $\ell_q$-(quasi-)norms
\begin{equation}
\left(\sum_i |a_i|^q\right)^{1/q} \leq \left(\sum_i |a_i|^p\right)^{1/p},
\text{ for } 0<p\leq q<\infty
\label{eq:lplqnorms}
\end{equation}
then provides
$$\sum_{I\in\mathcal{C}_{k+1}(K)}\|(\Sigma_{k+1}(s,K)-\Sigma_k(s,K))\BBone_I\|^q_p 
\leq 2^{q-1}\kappa^q(p) \left(e^q_{k+1}(s,K) + 2^{d(1+\underline{\bb \sigma}/H(\bb \sigma))} e^q_k(s,K)\right).$$
Since $(e_k(s,K))_k\in\N$ is a decreasing sequence, by setting $\tau=H(\bb \sigma)/d+1/q-1/p$  and combining Inequality~\eqref{eq:decompLq} with the above inequality, we obtain
\begin{align*}
\|\Sigma_{k+1}(s,K)-\Sigma_k(s,K)\|_q
&\leq C(d,\bb r,\bb \sigma,p,q) \left(2^{k\underline{\bb \sigma}}e_k(s,K)\right) 2^{-kd\tau\underline{\bb \sigma} /H(\sigma)} \\
&\leq C(d,\bb r,\bb \sigma,p,q) R 2^{-kd\tau\underline{\bb \sigma} /H(\sigma)}. 
\end{align*}
We can prove in the same way that such an upper-bound still holds for $q=\infty$. Since $\tau>0$, for all $0<p<q\leq\infty$, $(\Sigma_k(s,K))_{k\geq j}$ also converges in $\L_q([0,1]^d)$ to $s\BBone_K$. In particular, we have thus proved that $s\in\L_q([0,1]^d).$

From the definition of $\mathcal{E}_{\bb r}(s,K)_q$ and the triangle inequality, it follows that 
\begin{align}
\mathcal{E}_{\bb r}(s,K)_q
&\leq \|(s-P_K(s))\BBone_K\|_q \notag\\
&\leq \sum_{k\geq j}\|\Sigma_{k+1}(s,K)-\Sigma_k(s,K)\|_q\notag\\
&\leq C(d,\bb r,\bb \sigma,p,q) \sum_{k\geq j} 2^{-kd\tau\underline{\bb \sigma} /H(\sigma)} 2^{k\underline{\bb \sigma}} e_k(s,K).
\label{eq:changenormrectPK}
\end{align}
We have thus proved~\eqref{eq:changenormrect}, and the above inequality for $K=[0,1]^d$ combined with Markov Inequality~\eqref{eq:Markov} also provides 
\begin{align*}
\|s\|_q 
&\leq \|P_{[0,1]^d}(s)\|_q + \|s-P_{[0,1]^d}(s)\|_q  \\
&\leq C(d,\bb r,\bb \sigma,p,q) (\|P_{[0,1]^d}(s)\|_p+ R) \\
&\leq C(d,\bb r,\bb \sigma,p,q) (\|s\|_p + \mathcal{E}_{\bb r}(s,[0,1]^d)_p+R)\\
&\leq C(d,\bb r,\bb \sigma,p,q) (\|s\|_p +R).
\end{align*}
$\blacksquare$

\subsection{Proofs of Theorems~\ref{ApproxAnisotropicNC}  and~\ref{ApproxAnisotropicC}}

A first approximation result for the algorithm decsribed in Section~\ref{sec:approximation} can be stated as follows.
\begin{prop}\label{firstapprox}
Let $k\in\N$, $R>0$, $\bb \sigma=(\sigma_1,\ldots,\sigma_d) \in \prod_{l=1}^d (0,r_l+1)$, $0<p <\infty$, $1\leq q\leq \infty$ and $s\in\L_q([0,1]^d)$. 
Assume that 
$$H(\bb \sigma)/d>(1/p-1/q)_+$$
and that
\begin{equation}
\sup_{j\in\N}2^{jd\left(\underline{\bb \sigma}/H(\bb \sigma)\right) \left(H(\bb \sigma)/d-(1/p-1/q)_+\right)}\left(\sum_{K\in\mathcal{D}_j^{\bb \sigma}} \mathcal{E}_{\bb r}^p(s,K)_q\right)^{1/p} \leq R.
\label{eq:polyapprox}
\end{equation} 
Then, there exists some partition $m$ of $[0,1]^d$ that only contains dyadic rectangles from $\mathcal{D}^{\bb \sigma}$ and $s_\mr\in S_\mr$ such that 
$$|m|\leq C_1(d,\bb \sigma,p) 2^{kd}$$
and 
$$\|s-s_\mr\|_q \leq C_2(d,\bb \sigma,p,q) R 2^{-kH(\bb \sigma)}.$$
Besides, if for some $J\in\N$, $s$ is polynomial with coordinate degree $\leq r$ over each rectangle of $\mathcal{D}_J^{\bb \sigma}$, then $m$ only contains dyadic rectangles from $\cup_{j=0}^J\mathcal{D}_j^{\bb \sigma}$.
\end{prop}

\noindent
\textbf{\textit{Proof: }} For $k=0$, we can just choose $m$ as the trivial partition of $[0,1]^d$ and $s_\mr$ as the polynomial of best $\L_q$-approximation over $[0,1]^d$ in $\mathscr P_{\bb r}$. Indeed, we then have 
$$ \|s-s_\mr\|_q =\mathcal{E}_{\bb r}(s,[0,1]^d)_q \leq R,$$
where the last inequality follows from~\eqref{eq:polyapprox}. Let us now fix $k\geq 1$, set
$$\tau=H(\bb \sigma)/d-(1/p-1/q)_+ \quad\text{ and }\quad \lambda=2^{\left(1+(1+\tau p)\underline{\bb \sigma}/H(\bb \sigma)\right)d/p},$$
and choose 
$$\epsilon=\lambda R 2^{-kd\left(\tau+1/p\right)}.$$
If $\mathcal{I}(s,\epsilon)$ is trivial, then the upper-bound~\eqref{eq:simpleapproxfinite} provides
$$\|s-A(s,\epsilon)\|_q \leq \epsilon \leq \lambda R 2^{-kH(\bb \sigma)}.$$
Let us now assume that $\mathcal{I}(s,\epsilon)$ is not trivial and fix $j\geq 1$ such that $\mathcal{I}(s,\epsilon)\cap \mathcal{D}^{\bb \sigma}_j$ is not empty. If $K\in\mathcal{I}(s,\epsilon)\cap \mathcal{D}^{\bb \sigma}_j$, then $K$ is a child of a dyadic rectangle $\tilde K\in\mathcal{D}^{\bb \sigma}_{j-1}$ such that
$$\epsilon \leq \mathcal{E}_{\bb r}(s,\tilde K)_q,$$
hence 
$$\epsilon^p \leq 2^{-(j-1)dp\tau \underline{\bb \sigma} /H(\bb \sigma)} 
2^{(j-1)dp\tau \underline{\bb \sigma} /H(\bb \sigma)} \mathcal{E}_{\bb r}^p(s,\tilde K)_q.$$
By grouping all the rectangles $K\in\mathcal{I}(s,\epsilon)\cap\mathcal{D}_j^{\bb \sigma}$ having the same parent in $\mathcal{D}_{j-1}^{\bb \sigma}$, and taking into account Remark~\eqref{eq:children}, we obtain
\begin{equation}
|\mathcal{I}(s,\epsilon)\cap\mathcal{D}_j^{\bb \sigma}|\epsilon^p
\leq 2^{d(1+(1+p\tau)\underline{\bb \sigma}/H(\bb \sigma))} 2^{-jdp\tau\underline{\bb \sigma}/H(\bb \sigma)} R^p.
\notag
\end{equation}
Replacing $\epsilon$ by its value, we deduce that 
\begin{equation}
|\mathcal{I}(s,\epsilon)\cap\mathcal{D}_j^{\bb \sigma}|
\leq 2^{kd(1+p\tau)}2^{-jdp\tau\underline{\bb \sigma}/H(\bb \sigma)} .
\label{eq:boundNj}
\end{equation}
Besides, for all $j\geq 1$, 
\begin{equation}
|\mathcal{I}(s,\epsilon)\cap\mathcal{D}_j^{\bb \sigma}|
\leq |\mathcal{D}_j^{\bb \sigma}|
\leq 2^{jd\underline{\bb \sigma}/H(\bb \sigma)}.
\notag
\end{equation}
Let us denote by $J$ the greatest integer $j\geq 1$ such that
$$2^{jd\underline{\bb \sigma}/H(\bb \sigma)}\leq 2^{kd(1+p\tau)}2^{-jdp\tau\underline{\bb \sigma}/H(\bb \sigma)} ,$$
\textit{i.e.} such that
$$2^{jd\underline{\bb \sigma}/H(\bb \sigma)}\leq 2^{kd}.$$
Since $\underline{\bb \sigma}/H(\bb \sigma) \leq 1$, the last inequality is satisfied by $k\geq 1$ for instance, so that $J$ is well-defined. Besides, $J$ is characterized by 
$$2^{Jd\underline{\bb \sigma}/H(\bb \sigma)}\leq 2^{kd} < 2^{(J+1)d\underline{\bb \sigma}/H(\bb \sigma)}.$$
Therefore,
\begin{align}
|\mathcal{I}(s,\epsilon)|
&=\sum_{j\geq 1}|\mathcal{I}(s,\epsilon)\cap\mathcal{D}_j^{\bb \sigma}|\notag\\
&\leq \sum_{j=1}^J 2^{jd\underline{\bb \sigma}/H(\bb \sigma)} + 
2^{kd(1+p\tau)} \sum_{j\geq J+1} 2^{-jdp\tau\underline{\bb \sigma}/H(\bb \sigma)}\notag\\
&\leq C_1(d,\bb \sigma,p) 2^{kd}
\notag
\end{align}
where 
$$C_1(d,\bb \sigma,p) = \frac{2^{d\underline{\bb \sigma}/H(\bb \sigma)}}{2^{d\underline{\bb \sigma}/H(\bb \sigma)}-1} + \frac{1}{1-2^{-dp\tau \underline{\bb \sigma}/H(\bb \sigma)}}.$$
Moreover, we deduce from~\eqref{eq:simpleapproxfinite} that, if $1\leq q <\infty$, then
\begin{equation}
\|s-A(s,\epsilon)\|_q
\leq |\mathcal{I}(s,\epsilon)|^{1/q} \epsilon
\leq C^{1/q}_1(d,\bb \sigma,p)  R 2^{-kH(\bb \sigma)},
\notag
\end{equation}
and we deduce from~\eqref{eq:simpleapproxinfinite} that, if $q=\infty$, then
\begin{equation}
\|s-A(s,\epsilon)\|_\infty
<\epsilon
\leq \lambda  R 2^{-kH(\bb \sigma)}.
\notag
\end{equation}
So Proposition~\ref{firstapprox} is satisfied for 
$$C_2(d,\bb \sigma,p,q)=
\begin{cases}
C^{1/q}_1(d,\bb \sigma,p) \lambda &\text{ if } 1\leq q <\infty\\
\lambda &\text{ if } q=\infty,\\
\end{cases}$$
$m=\mathcal{I}(s,\epsilon)$ and $s_\mr=A(s,\epsilon)$. The last assertion in Proposition~\ref{firstapprox} is a straightforward consequence of the approximation algorithm. $\blacksquare$ 

\medskip
The following lemma allows to link Assumption~\eqref{eq:polyapprox} with the (quasi-)semi-norm $N_{\bb r,\bb \sigma,p,p'}.$ 
\begin{lemm}\label{linknorm}
Let $R>0$, $0<p<\infty$, $1\leq q \leq \infty$ and $\bb \sigma=(\sigma_1,\ldots,\sigma_d) \in \prod_{l=1}^d (0,r_l+1)$ such that $H(\bb \sigma)/d>(1/p-1/q)_+.$ 
Assume that $s\in\mathcal{S}(\bb r,\bb \sigma,p,p',R)$, where $p'=\infty$ if $0<p\leq 1$ or $p\geq q$  and $p'=p$ if $1<p<q$, then 
\begin{equation}
\sup_{j\in\N}2^{jd\left(\underline{\bb \sigma}/H(\bb \sigma)\right) \left(H(\bb \sigma)/d-(1/p-1/q)_+\right)}\left(\sum_{K\in\mathcal{D}_j^{\bb \sigma}} \mathcal{E}_{\bb r}^p(s,K)_q\right)^{1/p}
\leq C(d,\bb r,\bb \sigma,p,q)R.
\label{eq:tBesov}
\end{equation}
\end{lemm}

\noindent
\textbf{\emph{Proof: }} If $p\geq q$, then the left-hand side of Inequality~\eqref{eq:tBesov} is upper-bounded by  
\begin{equation*}
\sup_{j\in\N}2^{j\underline{\bb \sigma}}\left(\sum_{K\in\mathcal{D}_j^{\bb \sigma}} \mathcal{E}_{\bb r}^p(s,K)_p\right)^{1/p} 
= \sup_{j\in\N}2^{j\underline{\bb \sigma}} e_{\bb r,\bb \sigma,p,j}(s) \leq R.
\end{equation*}
Let us now assume that $p<q$ and set $\tau=H(\bb \sigma)/d+1/q-1/p$. From Inequality~\eqref{eq:changenormrect} in Proposition~\ref{changenormaniso}, we deduce that
\begin{multline}
2^{jd\tau\underline{\bb \sigma}/H(\bb \sigma)}\left(\sum_{K\in\mathcal{D}_j^{\bb \sigma}} \mathcal{E}_{\bb r}^p(s,K)_q\right)^{1/p}\\
\leq C(d,r,\bb \sigma,p,q) 2^{jd\tau\underline{\bb \sigma}/H(\bb \sigma)} 
\left(\sum_{K\in\mathcal{D}_j^{\bb \sigma}} \left(\sum_{k\geq j} 2^{-kd\tau\underline{\bb \sigma}/H(\bb \sigma)} 2^{k\underline{\bb \sigma}}e_k(s,K)\right)^p\right)^{1/p}.
\label{eq:towardspolyapprox}
\end{multline}
If $0<p\leq 1$, then the classical inequality between $\ell_p$ and $\ell_1$-(quasi-)norms recalled in~\eqref{eq:lplqnorms} leads to 
\begin{align*}
&2^{jd\tau\underline{\bb \sigma}/H(\bb \sigma)}\left(\sum_{K\in\mathcal{D}_j^{\bb \sigma}} \mathcal{E}_{\bb r}^p(s,K)_q\right)^{1/p}\\
&\leq C(d,\bb r,\bb \sigma,p,q) 2^{jd\tau\underline{\bb \sigma}/H(\bb \sigma)} 
\left(\sum_{K\in\mathcal{D}_j^{\bb \sigma}} 
\sum_{k\geq j} 2^{-kpd\tau\underline{\bb \sigma}/H(\bb \sigma)} 2^{kp\underline{\bb \sigma}}e^p_k(s,K)\right)^{1/p}\\
&\leq C(d,\bb r,\bb \sigma,p,q) 2^{jd\tau\underline{\bb \sigma}/H(\bb \sigma)} 
\left(\sum_{k\geq j} 2^{-kpd\tau\underline{\bb \sigma}/H(\bb \sigma)} 2^{kp\underline{\bb \sigma}}e^p_k(s,[0,1]^d)\right)^{1/p}\\
&\leq C(d,\bb r,\bb \sigma,p,q) \sup_{k\geq j}\left(2^{k\underline{\bb \sigma}}e_k(s,[0,1]^d)\right)
2^{jd\tau\underline{\bb \sigma}/H(\bb \sigma)} \left(\sum_{k\geq j} 2^{-kpd\tau\underline{\bb \sigma}/H(\bb \sigma)}\right)^{1/p}
\end{align*}
hence Inequality~\eqref{eq:tBesov}.
If $1<p<\infty$, then there exists $1<p^\star<\infty$ such that $1/p+1/p^\star=1$, so we obtain by applying H\"older inequality to~\eqref{eq:towardspolyapprox} that
\begin{align*}
&2^{jd\tau\underline{\bb \sigma}/H(\bb \sigma)}\left(\sum_{K\in\mathcal{D}_j^{\bb \sigma}} \mathcal{E}_{\bb r}^p(s,K)_q\right)^{1/p}\\
&\leq C(d,\bb r,\bb \sigma,p,q) 2^{jd\tau\underline{\bb \sigma}/H(\bb \sigma)} 
\left(\sum_{K\in\mathcal{D}_j^{\bb \sigma}} \left(\sum_{k\geq j} 2^{-p^\star kd\tau\underline{\bb \sigma}/H(\bb \sigma)} \right)^{p/p^\star}  \left(\sum_{k\geq j} 2^{kp\underline{\bb \sigma}}e^p_k(s,K)\right)\right)^{1/p}\\
&\leq C(d,\bb r,\bb \sigma,p,q) 
 \left(\sum_{k\geq j} 2^{kp\underline{\bb \sigma}}\sum_{K\in\mathcal{D}_j^{\bb \sigma}}e^p_k(s,K)\right)^{1/p}\\
&\leq C(d,\bb r,\bb \sigma,p,q) 
 \left(\sum_{k\geq j} 2^{kp\underline{\bb \sigma}}e^p_k(s,[0,1]^d)\right)^{1/p}\\
\end{align*}
hence Inequality~\eqref{eq:tBesov}. $\blacksquare$

Last, Lemma~\ref{linrate} provides an upper-bound for the linear approximation error of $\mathcal{S}(\bb r,\bb \sigma,p,p',R)$ by $\Pi^{\bb r,\bb \sigma}_J$ in the $\L_q$-norm.

\begin{lemm}\label{linrate}
Let $R>0$, $0<p<\infty$, $1\leq q\leq \infty$, $\bb \sigma=(\sigma_1,\ldots,\sigma_d) \in \prod_{l=1}^d (0,r_l+1)$ such that $H(\bb \sigma)/d>(1/p-1/q)_+$, and $\kappa(p)=2^{1/p}$ if $0<p\leq 1$,  and $1$ otherwise.
Assume that $s\in\mathcal{S}(\bb r,\bb \sigma,p,p',R)$ where $p'=\infty$ if $0<p\leq 1$ or $p\geq q$, and $p'=p$ if $1<p<q$. Then, for all $J\in\N$, there exists a function $s_J\in\Pi^{\bb r,\bb \sigma}_J$ such that $s_J\in\mathcal{S}(\bb r,\bb \sigma,p,p',2\kappa(p)R)$ and
\begin{equation}
\|s-s_J\|_q \leq C(d,\bb r,\bb \sigma,p,q) 2^{-Jd(H(\bb \sigma)/d-(1/p-1/q)_+) \underline{\bb \sigma}/H(\bb \sigma)}R.
\label{eq:linapprox}
\end{equation} 
\end{lemm}

\noindent
\textbf{\emph{Proof: }} For all $K\in\mathcal{D}_J^{\bb \sigma}$, we denote by $P_K(s)$ a polynomial function on $K$ with degree $\leq r_l$ in the $l$-th direction such that
$$\|(s-P_K(s))\BBone_I\|_p = \mathcal{E}_{\bb r}(s,K)_p,$$
and we set 
$$s_J=\sum_{K\in\mathcal{D}^{\bb \sigma}_J}P_K(s)\BBone_K.$$
In order to alleviate the notation, we simply write $e_k(s)$ instead of $e_{\bb r,\bb \sigma,p,k}(s)$, and  $e_k(s,K)$ instead of $e_{\bb r,\bb \sigma,p,k}(s,K)$, as in the proof of Proposition~\ref{changenormaniso}. 

Since $s_J\in\Pi^{\bb r, \bb \sigma}_J$, $e_k(s_J)=0$ for $k\geq J$. If $k<J$, then the (quasi-)triangle inequality, the definition of $s_J$ and the inclusion $\Pi^{\bb r,\bb \sigma}_k \subset \Pi^{\bb r,\bb \sigma}_J$ provide successively 
\begin{align}
e_k(s_J)
&\leq \kappa(p) \left(\|s-s_J\|_p+e_k(s)\right)\notag\\
&\leq \kappa(p) \left(e_J(s)+e_k(s)\right)\notag\\
&\leq 2\kappa(p) e_k(s).
\notag
\end{align}
Therefore, $N_{\bb r,\bb \sigma,p,p'}(s_J) \leq 2\kappa(p) N_{\bb r,\bb \sigma,p,p'}(s)$, so that $s_J\in\mathcal{S}(\bb r,\bb \sigma,p,p',2\kappa(p)R).$

If $p\geq q$, then $$\|s-s_J\|_q \leq  \|s-s_J\|_p
= \left(\sum_{K\in\mathcal{D}^{\bb \sigma}_J}\mathcal{E}^p_{\bb r}(s,K)_p\right)^{1/p}
= e_{\bb r,\bb \sigma,p,J} (s)
\leq 2^{-J\underline{\bb \sigma}} R.$$
If $p<q<\infty$, then we deduce from Inequality~\eqref{eq:changenormrectPK} in the proof of Proposition~\ref{changenormaniso} and from Inequality~\eqref{eq:lplqnorms} between $\ell_p$ and $\ell_q$-(quasi-)norms that
\begin{align*}
\|s-s_J\|_q 
&= \left(\sum_{K\in\mathcal{D}^{\bb \sigma}_J} \|(s-P_K(s))\BBone_K\|^q_q\right)^{1/q}\\
&\leq C(d,\bb r,\bb \sigma,p,q) \left( \sum_{K\in\mathcal{D}^{\bb \sigma}_J}\left(\sum_{k\geq J} 2^{-kd(H(\sigma)/d+1/q-1/p)\underline{\bb \sigma}/H(\bb \sigma)}2^{k\underline{\bb \sigma}} e_k(s,K)\right)^p\right)^{1/p}.
\end{align*}
We then obtain Inequality~\eqref{eq:linapprox} either thanks to the inequality between $\ell_1$ and $\ell_p$-(quasi-)norms in case $0\leq p\leq 1$, or thanks to H\"older Inequality otherwise. Last, if $q=\infty$, then we still deduce from Inequality~\eqref{eq:changenormrectPK} that  
\begin{align*}
\|s-s_J\|_\infty
&= \max_{K\in\mathcal{D}^{\bb \sigma}_J} \|(s-P_K(s))\BBone_K\|_\infty\\
&\leq C(d,\bb r,\bb \sigma,p,q)  \max_{K\in\mathcal{D}^{\bb \sigma}_J} \left(\sum_{k\geq J} 2^{-kd(H(\sigma)/d+1/q-1/p)\underline{\bb \sigma}/H(\bb \sigma)}2^{k\underline{\bb \sigma}} e_k(s,K)\right)\\
&\leq C(d,\bb r,\bb \sigma,p,q)  2^{-Jd(H(\sigma)/d+1/q-1/p)\underline{\bb \sigma}/H(\bb \sigma)} R.
\end{align*}
$\blacksquare$

Theorem~\ref{ApproxAnisotropicNC} is then a straightforward consequence of Proposition~\ref{firstapprox} and Lemma~\ref{linknorm}. To prove Theorem~\ref{ApproxAnisotropicC}, for each $J\in\N$, we just have to apply Proposition~\ref{firstapprox} and Lemma~\ref{linknorm} to the function $s_J$ given by Lemma~\ref{linrate} and use the triangle inequality 
\begin{equation*}
\inf_{t\in S_\mr} \|s-t\|_q \leq \|s-s_J\|_q + \inf_{t\in S_\mr}\|s_J-t\|_q
\end{equation*}
where $m$ can be any partition of $[0,1]^d$ into dyadic rectangles.


\section{Proof of Theorem~\ref{selection}}\label{sec:proofselect}

In the following proof, we denote by $(w_\mrho)_{\mrho\in\Mdeg}$ a family of nonnegative reals and set $\Sigma=\sum_{\mrho\in\Mdeg} \exp(-w_\mrho)$. We fix $\mrho\in\Mdeg$ as well as some positive reals $\zeta$, $\theta_1,\ldots,\theta_8$ such that $2\theta_1(1+\theta_2)<1$ and $\theta_8<1$.

From the definition of $\tilde s= \hat s _\hatmrho$, it follows  that 
\begin{equation}
\gamma(\tilde s) + \pen\hatmrho \leq \gamma(\hat s_\mrho) + \pen\mrho.
\label{eq:key}
\end{equation}
For all $t,u\in\L_2([0,1]^d)$,
\begin{equation}
\gamma(t)-\gamma(u)=\|s-t\|_2^2 - \|s-u\|_2^2 - 2 \nu(t-u),
\label{eq:decompgamma}
\end{equation}
where 
$$\nu(t)=\frac{1}{n} \sum_{i=1}^n \left(t(Y_i) - \langle t,s \rangle \right).$$
Besides, for all $\primemrho\in\Mdeg,$ setting 
$$\chi\primemrho=\|s_\primemrho -\hat s_\primemrho\|_2,$$
we obtain by developing $s_\primemrho$ and $\hat s_\primemrho$ in the orthonormal basis $(\Phi_{K,\bb k})_{K\in m', \bb k\in\Lambda(\bb{\rho}'_K)}$ and using the linearity of $\nu$ 
\begin{equation}
\chi^2\primemrho
= \sum_{K\in m'} \sum_{\bb k\in\Lambda(\bb{\rho}'_K)}  \nu^2(\Phi_{K,\bb k})
= \nu\left(\hat s_\primemrho - s_\primemrho\right).
\label{eq:chinu}
\end{equation}
From Equalities~\eqref{eq:decompgamma},~\eqref{eq:chinu}, Pythagoras' Equality and the linearity of $\nu$, we deduce 
\begin{equation*}
\gamma(\tilde s)-\gamma(\hat s_\mrho)
=\|s-\tilde s\|_2^2 - \|s- s_\mrho\|_2^2 + \chi^2\mrho-2\chi^2\hatmrho-2\nu\left(s_\hatmrho - s_\mrho\right),
\end{equation*}
which, combined with Inequality~\eqref{eq:key}, leads to
\begin{align}
\|s-\tilde s\|_2^2 &\leq \|s- s_\mrho\|_2^2 + \pen\mrho-\chi^2\mrho\notag\\
& +2\chi^2\hatmrho+2\nu\left(s_\hatmrho - s_\mrho\right) -\pen\hatmrho.
\label{eq:step1}
\end{align}

We shall now provide an upper-bound for the term $2\nu\left(s_\hatmrho - s_\mrho\right)$ on an event with great probability. From Bernstein's Inequality, as stated for instance in~\cite{Massart} (Section 2.2.3), for all bounded function $t:[0,1]^d \rightarrow \R$ and all $x>0$, 
\begin{equation}
\P_s\left(\nu(t) \geq \sqrt{2\E_s[t^2(Y_1)] \frac{x}{n}} + \frac{\|t\|_\infty}{3}\frac{x}{n} \right) \leq \exp(-x).
\label{eq:Bernsteinonesided}
\end{equation}
Let us fix $\primemrho\in\Mdeg$ and apply Bernstein's Inequality to $t=s_\primemrho - s_\mrho$. Since $\Phi_{K,\bb k}$ has support $K$,
$$\|s_\primemrho\|_\infty = \max_{K\in m'} \left\|\sum_{\bb k\in\Lambda(\bb{\rho}'_K)} \langle s, \Phi_{K,\bb k}  \rangle \Phi_{K,\bb k} \right\|_\infty \leq \Mone$$
where
$$\Mone=\max_{K\in\mathcal D_\star} \sum_{\bb k\in \Lambda(\bb{r_\star})} \sqrt{\frac{\pi(\bb k)}{\lambda_d(K)}} |\langle s, \Phi_{K,\bb k}  \rangle|,$$
so 
$$\|s_\primemrho-s_\mrho\|_\infty  \leq 2\Mone.$$
Since $S_\mrho$ and $S_\primemrho$ are both subspaces of $S_\starmr$, 
\begin{align*}
\E_s\left[\left(s_\primemrho - s_\mrho\right)^2(Y_1)\right] 
&= \int_{[0,1]^d} s_\starmr\left(s_\primemrho - s_\mrho\right)^2 \d\lambda_d \\
&\leq \Mone \|s_\primemrho - s_\mrho\|_2^2.
\end{align*}
From~\eqref{eq:Bernsteinonesided}, there exists a set $\Omega(m,\bb\rho,m',\bb{\rho'},\zeta)$ such that $\P_s\left(\Omega(m,\bb\rho,m',\bb{\rho'},\zeta)\right)\geq 1-\exp(-(w_\primemrho+\zeta))$ and over which
$$\nu\left(s_\primemrho - s_\mrho\right) 
\leq \sqrt{2\Mone \|s_\primemrho - s_\mrho\|_2^2\frac{w_\primemrho+\zeta}{n}} + \frac{2}{3}\Mone\frac{w_\primemrho+\zeta}{n}.$$  
We recall that, for all $a,b\geq 0$ and $\theta >0$, 
\begin{equation}
2ab \leq \theta a^2 + \theta^{-1} b^2.
\label{eq:ab}
\end{equation}
Thus, on $\Omega(m,\bb\rho,m',\bb{\rho'},\zeta)$, we have
\begin{equation*}
\nu\left(s_\primemrho - s_\mrho\right) 
\leq \theta_1\|s_\primemrho - s_\mrho\|_2^2 + \left(2/3+\theta_1^{-1}\right)\Mone\frac{w_\primemrho+\zeta}{n}.
\end{equation*}
Besides, using the triangle inequality,~\eqref{eq:ab}, and Pythagoras' Equality, we obtain 
\begin{align*}
\|s_\primemrho - s_\mrho\|_2^2 
&\leq \left(\|s-s_\primemrho\|_2+\|s - s_\mrho\|_2\right)^2\\
&\leq (1+\theta_2)\|s-s_\primemrho\|_2^2+\left(1+\theta_2^{-1} \right)\|s - s_\mrho\|_2^2\\
&\leq (1+\theta_2)\|s-\hat s_\primemrho\|_2^2- (1+\theta_2)\chi^2\primemrho+\left(1+\theta_2^{-1}\right)\|s - s_\mrho\|_2^2.
\end{align*}
Therefore, the set $\Omega_\mrho(\zeta)=\cap_{\primemrho \in \Mdeg} \Omega(m,\bb\rho,m',\bb{\rho'},\zeta)$ is an event with probability
\begin{equation}
\P_s\left(\Omega_\mrho(\zeta)\right) \geq 1-\exp(-\zeta) \Sigma
\label{eq:probanu}  
\end{equation}
over which
\begin{multline}
2\nu\left(s_\hatmrho - s_\mrho\right) 
\leq 2\theta_1 (1+\theta_2)\|s-\tilde s\|_2^2 + 2\theta_1\left(1+\theta_2^{-1} \right)\|s - s_\mrho\|_2^2
\\- 2\theta_1(1+\theta_2)\chi^2\hatmrho+2\left(2/3+\theta_1^{-1}\right)\Mone\frac{w_\hatmrho+ \zeta}{n}.
\label{eq:controlnu}  
\end{multline}

Let us now provide a concentration inequality for $\chi^2\hatmrho$. For that purpose, we first prove the following result.
\begin{prop}\label{concentrationchi}
Let $\primemrho\in\Mdeg,x>0$,
$$V_\primemrho
=\frac{1}{n}\sum_{K\in m'}\sum_{\bb k \in \Lambda({\bb\rho}'_K)} \Var_s\left(\Phi^2_{K,\bb k}(Y_1)\right)
=\E_s\left[\|\hat s_\primemrho - s_\mrho\|_2^2\right]$$
and
$$M_{2,\star}=\max_{K\in\mathcal D_\star,\bb k \in \Lambda (\bb r_\star)} \E_s[\Phi^2_{K,\bb k}(Y_1)].$$ There exist an event $\Omega_\star$ that does not depend on $\primemrho$ and an event  $\Omega_\primemrho(x)$ such that 
\begin{equation}
\P_s(\Omega_\star^c) \leq 2^{d+1}/(n^2\log(n)),
\label{eq:omegastar}
\end{equation}
\begin{equation*}
\P_s(\Omega^c_\primemrho(x)) \leq \exp(-x),
\label{eq:omegachimod}
\end{equation*}
and, on $\Omega_\primemrho(x)$,
\begin{multline}
\chi^2\primemrho \BBone_{\Omega_\star} \leq (1+\theta_3)(1+\theta_4)V_\primemrho \\
+ 4\left(1+\theta_4^{-1}\right)|\Lambda(\bb r_\star)|
\left(\left(4/3+\theta_3^{-1}\right)M_{2,\star}+(5/3)\left(1+3\theta_3^{-1}\right)\pi(\bb r_\star)\right)\frac{x}{n}.
\notag
\end{multline}
\end{prop}
\textit{\textbf{Proof: }}Let us fix $x>0$, and set, for all $K\in \mathcal D_\star$ and $\bb k\in \Lambda(\bb r_\star)$, 
$$\sigma^2_{K,k}=\Var_s\left(\Phi_{K, \bb k} (Y_1)\right),\quad \varepsilon_{K,\bb k}=\sqrt{6\sigma^2_{K,k}} +2 \sqrt{\pi(\bb k)},$$
$$\Omega_\star=\bigcap_{K\in \mathcal D_\star} \bigcap_{\bb k \in \Lambda(\bb r_\star)} \left\{|\nu(\Phi_{K,\bb k})|< \varepsilon_{K,\bb k} \sqrt{\lambda_d(K)}\right\}.$$ 
From Bernstein's Inequality (see for instance~\cite{Massart}, Section 2.2.3), for all $K\in \mathcal D_\star$, $\bb k\in \Lambda(\bb r_\star)$ and $x>0$,
$$\P_s\left(|\nu(\Phi_{K,\bb k})|\geq \sqrt{2\sigma^2_{K,\bb k} \frac{x}{n}}+ \frac{2\sqrt{\pi(\bb k)}}{3\sqrt{\lambda_d(K)}}\frac{x}{n}\right) \leq 2\exp(-x),$$
so  
$$\P_s\left(|\nu(\Phi_{K,\bb k})|\geq \varepsilon_{K,\bb k}\sqrt{\lambda_d(K)}\right)\leq 2\exp(-3 n \lambda_d(K))\leq 2\exp(-3n 2^{-dJ_\star}).$$
Besides, there are $2^{J_\star+1}-1$ dyadic intervals of $[0,1]$ with length $\geq 2^{-J_\star}$, so $|\mathcal D_\star|\leq 2^{d(1+J_\star)}.$ And we assume that $2^{dJ_\star}\leq n/\log(n|\Lambda(\bb r_\star)|)$, hence the upper-bound for $\P_s(\Omega_\star^c).$

Let us also fix $\primemrho\in\Mdeg$, set
$$v_\primemrho=\max_{K\in m'}\sum_{k\in\Lambda(\bb{\rho}'_K)} \varepsilon_{K,\bb k} \sqrt{\pi(\bb k)}
\quad\text{ and }\quad
b_\primemrho(x)=\sqrt{\frac{nv_\primemrho}{2(1/3+\theta_3^{-1})x}},$$
choose $\mathscr T_\primemrho$ a countable and dense subset of $\mathcal T_\primemrho=\left\{t\in S_\primemrho /\|t\|_2=1, \|t\|_\infty\leq b_\primemrho(x)\right\},$
and define
$$Z\primemrho=\sup_{t\in\mathcal T_\primemrho} \nu(t)=\sup_{t\in\mathscr T_\primemrho} \nu(t).$$ 
Since $\Phi_{K,\bb k}$ has support $K$, for all $t\in S_\primemrho,$ 
\begin{align*}
\E_s\left[t^2(Y_1)\right]
&=\E_s\left[\sum_{K\in m'}\left(\sum_{\bb k \in \Lambda(\bb{\rho'}_K)} \langle t,\Phi_{K,k}\rangle \Phi_{K,k} \right)^2\right]\\
&\leq\sum_{K\in m'}|\Lambda(\bb{\rho}'_K)|\sum_{\bb k \in \Lambda(\bb{\rho}'_K)} \langle t,\Phi_{K,k}\rangle^2 \E_s\left[\Phi^2_{K,k}(Y_1)\right]\\
&\leq|\Lambda(\bb r_\star)| M_{2,\star}\|t\|_2^2.
\end{align*}
So Talagrand's Inequality, as stated for instance in~\cite{Massart} (Chapter 5, Inequality (5.50)), ensures that there exists an event $\Omega_\primemrho(x)$ such that $\P_s(\Omega_\primemrho(x))\geq 1-\exp(-x)$ and over which 
$$Z\primemrho \leq 
(1+\theta_3)\E_s\left[Z\primemrho\right] +\sqrt{2|\Lambda(\bb r_\star)|M_{2,\star}\frac{x}{n}} +\sqrt{2(1/3+\theta_3^{-1})v_\primemrho\frac{x}{n}}.$$
Since $\nu$ is linear, we deduce from Cauchy-Scwharz Inequality and its equality case that 
$$\chi\primemrho=\sup_{t\in S_\primemrho, \|t\|_2=1} \nu(t) = \nu(t^\bullet_\primemrho)$$
where 
$$t^\bullet_\primemrho=\sum_{K\in m'}\sum_{\bb k \in \Lambda(\bb{\rho}'_K)} \frac{\nu(\Phi_{K,\bb k}) }{\chi\primemrho}\Phi_{K,\bb k}.$$
Therefore, $$\E_s\left[Z\primemrho\right] \leq \E_s\left[\chi\primemrho\right] \leq \sqrt{\E_s\left[\chi^2\primemrho\right]}=\sqrt{V_\primemrho}.$$
Moreover, on the set $\Omega_\primemrho(x)\cap \Omega_\star$,
either $\chi\primemrho \geq \sqrt{2(1/3+\theta_3^{-1})v_\primemrho x/n}$, in which case $t^\bullet_\primemrho\in \mathcal T_\primemrho$, so that 
\begin{align*}
\chi\primemrho
&=Z\primemrho\\
&\leq 
(1+\theta_3)\sqrt{V_\primemrho} +\sqrt{2|\Lambda(\bb r_\star)|M_{2,\star}\frac{x}{n}} +\sqrt{2(1/3+\theta_3^{-1})v_\primemrho\frac{x}{n}}, 
\end{align*}
or $\chi\primemrho < \sqrt{2(1/3+\theta_3^{-1})v_\primemrho x/n}$, and the above inequality is still satisfied. Applying Inequality~\eqref{eq:ab} with $\theta=1$, we get 
$$v_\primemrho
\leq \max_{K\in m'}\sum_{\bb k \in \Lambda(\bb{\rho}'_K)} \left(\sigma^2_{K,\bb k}+5 \pi(\bb k)\right)\leq |\Lambda(\bb r_\star)|\left(M_{2,\star}+5 \pi(\bb r_\star)\right).$$
Consequently, on $\Omega_\primemrho(x)$, 
\begin{multline}
\chi\primemrho\BBone_{\Omega_\star}
\leq (1+\theta_3)\sqrt{V_\primemrho} \\
+\left(\sqrt{M_{2,\star}}+\sqrt{(1/3+\theta_3^{-1})(M_{2,\star}+5 \pi(\bb r_\star)}\right)\sqrt{2|\Lambda(\bb r_\star)|\frac{x}{n}}.
\notag
\end{multline} 
Thus, applying twice Inequality~\eqref{eq:ab}, with $\theta=\theta_4$ and $\theta=1$, we get the concentration inequality for $\chi\primemrho$ stated in Proposition~\ref{concentrationchi}. $\blacksquare$
\\\noindent
From Proposition~\ref{concentrationchi}, we deduce that $\Omega_\chi(\zeta)=\cap_{\primemrho \in \Mdeg} \Omega_\primemrho(w_\primemrho+\zeta)$ is an event with probability 
\begin{equation*}
\P_s\left(\Omega_\chi(\zeta)\right)\geq 1-\exp(-\zeta)\Sigma
\end{equation*}
over which
\begin{multline}
\chi^2\hatmrho \BBone_{\Omega_\star} \leq (1+\theta_3)(1+\theta_4)V_\hatmrho \\
+ 4\left(1+\theta_4^{-1}\right)\left|\Lambda(\bb r_\star)\right|
\left(\left(4/3+\theta_3^{-1}\right)M_{2,\star}+(5/3)\left(1+3\theta_3^{-1}\right)\pi(\bb r_\star)\right)\frac{w_\hatmrho+\zeta}{n}.
\label{eq:controlchi}
\end{multline}

Our main task is then to estimate the unknown variance terms $V_\primemrho,\Mone,M_{2,\star}.$ Lemma 1 in~\cite{RBRTM} remains valid with the same constants even though the $Y_i$'s take values in $\R^d$ with $d\geq 1$. Let us set $\gamma=3+\log|\Lambda(\bb r_\star)|/\log(n)$. Since $|\Lambda(\bb r_\star)|\leq n^d$, $\gamma$ is bounded independently of $n$ ( $3\leq \gamma\leq 3+d)$). So, from the proof of Lemma 1 in~\cite{RBRTM}, for all $K\in\mathcal D_\star$ and $\bb k\in \Lambda(\bb r_\star)$, there exists an event $\Omega_{K,\bb k}$ such that $\P_s(\Omega^c_{K,\bb k})\leq C(\theta_5,d)/(n^3 |\Lambda(\bb r_\star)|)$ and over which
\begin{align*}
\Var_s\left(\Phi_{K,\bb k}(Y_1)\right)
& \leq (1+\theta_5)\left(\hat \sigma^2_{K,\bb k} + 2 \|\Phi_{K,\bb k}\|_\infty\sqrt{2\gamma\hat \sigma^2_{K,\bb k}\frac{\log(n)}{n}} + 8\gamma\|\Phi_{K,\bb k}\|^2_\infty \frac{\log(n)}{n}\right)\\
& \leq (1+\theta_5)\left(\hat \sigma^2_{K,\bb k} + 2 \sqrt{8\hat \sigma^2_{K,\bb k}\pi(\bb k)} + 32 \pi(\bb k)\right).
\end{align*}
Applying Inequality~\eqref{eq:ab} with $a=\hat \sigma_{K,\bb k}$, $b=\sqrt{8\pi(\bb k)}$ and $\theta=\theta_6$, we get, on  $\Omega_{K,\bb k}$,
$$\sigma^2_{K,\bb k} \leq (1+\theta_5) \left((1+\theta_6)\hat \sigma^2_{K,\bb k}+8(4+\theta_6^{-1})\pi(\bb k)\right).$$
For all $\primemrho\in\Mdeg$, let us introduce 
$$\widehat V_\primemrho (\theta_6)=\frac{1}{n}\sum_{K\in m'} \sum_{\bb k \in \Lambda(\bb{\rho}'_K)}\left((1+\theta_6)\hat \sigma^2_{K,\bb k}+8(4+\theta_6^{-1})\pi(\bb k)\right).$$ We have just proved that the set $\Omega_\sigma=\cap_{K\in\mathcal D_\star }\cap_{\bb k\in \Lambda(\bb r_\star)} \Omega_{K,\bb k}$ is an event with probability
\begin{equation}
\P_s\left(\Omega_\sigma\right) \geq 1-2^dC(\theta_5,d)/(n^2\log(n))
\label{eq:omegasigma}
\end{equation}
over which
\begin{equation}
V_\hatmrho \leq (1+\theta_5) \widehat V_\hatmrho (\theta_6).
\label{eq:controlV}
\end{equation}
Let us now fix $K\in \mathcal D_\star$ and $\bb k\in \Lambda(\bb r_\star)$. According to Bernstein's Inequality and Inequality~\eqref{eq:ab}, there exist events $\Omega_{K,\bb k}^1$ and $\Omega_{K,\bb k}^2$, each with $\P_s$-measure $\geq 1-2\exp(-3n\lambda_d(K))$, such that on $\Omega_{K,\bb k}^1$
\begin{align*}
\sqrt{\frac{\pi(\bb k)}{\lambda_d(K)}} \left|\frac{1}{n}\sum_{i=1}^n \Phi_{K,\bb k} (Y_i) - \E_s\left[\Phi_{K,\bb k} (Y_1)\right]\right| 
&\leq \sqrt{6\E_s\left[\Phi^2_{K,\bb k} (Y_1)\right]\pi(\bb k)} + \pi(\bb k)\\
&\leq \theta_7 \E_s\left[\Phi^2_{K,\bb k} (Y_1)\right] +(1+3\theta_7^{-1})\pi(\bb k),
\end{align*}
and on $\Omega_{K,\bb k}^2$
\begin{align*}
\left|\frac{1}{n}\sum_{i=1}^n \Phi^2_{K,\bb k} (Y_i) - \E_s\left[\Phi^2_{K,\bb k} (Y_1)\right]\right| 
&\leq \sqrt{6\|\Phi_{K,\bb k}\|^2_\infty \E_s\left[\Phi^2_{K,\bb k} (Y_1)\right]\lambda_d(\bb k)} + \|\Phi_{K,\bb k}\|^2_\infty \lambda_d(\bb k) \\
&\leq \theta_8  \E_s\left[\Phi^2_{K,\bb k} (Y_1)\right] +(1+3\theta_8^{-1})\pi(\bb k).
\end{align*}
We thus obtain that $\Omega_M=\cap_{K\in\mathcal D_\star}\cap_{\bb k\in \Lambda(\bb r_\star)}(\Omega_{K,\bb k}^1\cap\Omega_{K,\bb k}^2)$ is an event with probablity 
\begin{equation}
\P_s(\Omega_M)\geq 1-4\times 2^d/(n^2\log(n))
\label{eq:omegaM}
\end{equation}
over which
\begin{align}
\Mone&\leq \widehat \Mone + \theta_7(1-\theta_8)^{-1}|\Lambda(\bb r_\star)|\widehat M_{2,\star}+\left(\theta_7(1-\theta_8)^{-1}(1+3\theta_8^{-1})+(1+3\theta_7^{-1})\right)\pi(\bb r_\star)|\Lambda(\bb r_\star)|\notag\\
\widehat \Mone&\leq \Mone + \theta_7 |\Lambda(\bb r_\star)|M_{2,\star}+(1+3\theta_7^{-1})\pi(\bb r_\star)|\Lambda(\bb r_\star)|\notag\\
M_{2,\star} &\leq (1-\theta_8)^{-1}\widehat M_{2,\star}+(1+3\theta_8^{-1})(1-\theta_8)^{-1}\pi(\bb r_\star)\notag\\
\widehat M_{2,\star} &\leq (1+\theta_8) M_{2,\star}+(1+3\theta_8^{-1})\pi(\bb r_\star).
\label{eq:controlM}
\end{align}

Let us set $\Omega_\bullet=\Omega_\star\cap\Omega_\sigma\cap\Omega_M$ and 
\begin{align*}
C_0&= 1-2\theta_1(1+\theta_2)\\
C_1&= 1+2\theta_1(1+\theta_2^{-1})\\
C_2&=(1+C_0)(1+\theta_3)(1+\theta_4)(1+\theta_5)\\
C_3&=4(1+C_0)(4/3+\theta_3^{-1})(1+\theta_4^{-1})(1-\theta_8)^{-1}+C_7\theta_7(1-\theta_8)^{-1}\\
C_4&=3C_3+(20/3)(1+C_0)(1+3\theta_3^{-1})(1+\theta_4^{-1})+C_7\left(1+3\theta_7^{-1}+\theta_7(1+3\theta_8^{-1})(1-\theta_8)^{-1}\right)\\
C_5&=2(2/3+\theta_1^{-1})\\
C_6&=C_7+4(1+C_0)(4/3+\theta_3^{-1})(1+\theta_4^{-1})\\
C_7&=(20/3)(1+C_0)(1+3\theta_3^{-1})(1+\theta_4^{-1}).
\end{align*}
We choose $\pen$ such that, on $\Omega_\bullet$ and for all $\primemrho\in\Mdeg$,
\begin{equation*}
\pen\primemrho= C_2 \widehat V_\primemrho (\theta_6)
+\left(\left(C_3\widehat M_{2,\star}+ C_4\pi(\bb r_\star)\right) |\Lambda(\bb r_\star)|+C_5\widehat M_{1,\star}\right)\frac{w_\primemrho}{n}.
\end{equation*}
Thus, combining Inequalities~\eqref{eq:step1},~\eqref{eq:controlnu},~\eqref{eq:controlchi},~\eqref{eq:controlV},~\eqref{eq:controlM} with the upper-bounds 
$$M_{1,\star} \leq \pi(\bb r_\star) |\Lambda(\bb r_\star)|  \|s\|_\infty
\quad\text{and}\quad
M_{2,\star} \leq \pi(\bb r_\star) \|s\|_\infty,$$
we obtain, on $\Omega_m(\zeta)\cap\Omega_\chi(\zeta)\cap\Omega_\bullet$,
\begin{equation*}
C_0 \|s-\tilde s\|_2^2 
\leq C_1 \|s-s_\mrho\|_2^2 +\pen\mrho
+\left(C_6\|s\|_\infty+C_7\right)\pi(\bb r_\star)|\Lambda(\bb r_\star)|\frac{\zeta}{n}.
\end{equation*}
Setting 
\begin{align*}
C'_3&=C_3(1+\theta_8)+C_5(1+\theta_7)\\
C'_4&=C_3(1+3\theta_8^{-1})+C_5(1+3\theta_7^{-1})
\end{align*}
we deduce from~\eqref{eq:controlM} that, on $\Omega_\bullet$,
\begin{equation*}
\pen\mrho\leq C_2 \widehat V_\mrho(\theta_6)
+\left(C'_3\|s\|_\infty+ C'_4\right) \pi(\bb r_\star)|\Lambda(\bb r_\star)|\frac{w_\mrho}{n},
\end{equation*}
so that, on $\Omega_m(\zeta)\cap\Omega_\chi(\zeta),$
\begin{align*}
C_0 \|s-\tilde s\|_2^2\BBone_{\Omega_\bullet} 
&\leq C_1 \|s-s_\mrho\|_2^2 +C_2 \widehat V_\mrho(\theta_6)\\
&+\left(C'_3\|s\|_\infty+ C'_4 \right) \pi(\bb r_\star)|\Lambda(\bb r_\star)|\frac{w_\mrho}{n}\\
&+\left(C_6\|s\|_\infty+C_7\right)\pi(\bb r_\star)|\Lambda(\bb r_\star)|\frac{\zeta}{n}.
\end{align*}
Last, we recall that Fubini's Theorem yields, for all random variable $U$, 
$$\E[U]\leq \E[U_+] =\int_0^\infty\P(U_+>\zeta)\d \zeta= \int_0^\infty\P(U>\zeta)\d \zeta,$$
and we underline that 
$$\E_s\left[\widehat V_\mrho(\theta_6)\right]\leq (1+\theta_6) \E_s\left[\|\hat s_\mrho-s_\mrho\|_2^2\right] + 8(4+\theta_6^{-1})\pi(\bb r_\star)\frac{\dim(S_\mrho)}{n}.$$
Therefore, 
\begin{align}
C_0 \E_s\left[\|s-\tilde s\|_2^2\BBone_{\Omega_\bullet}\right] 
&\leq C_1 \|s-s_\mrho\|_2^2 
+(1+\theta_6)C_2\E_s\left[\|\hat s_\mrho-s_\mrho\|_2^2\right]\notag\\
&+ 8(4+\theta_6^{-1})C_2\pi(\bb r_\star)\frac{\dim(S_\mrho)}{n}\notag\\
&+\left(C'_3\|s\|_\infty+ C'_4\right) \pi(\bb r_\star)|\Lambda(\bb r_\star)|\frac{w_\mrho}{n}\notag\\
&+2\left(C_6\|s\|_\infty+C_7\right)\pi(\bb r_\star)|\Lambda(\bb r_\star)|\frac{\Sigma}{n}.
\label{eq:riskomega}
\end{align}

There remains to bound the risk of $\tilde s$ on $\Omega^c_\bullet$. According to~\eqref{eq:omegastar},~\eqref{eq:omegasigma} and~\eqref{eq:omegaM},
$$p_\bullet:=\P_s(\Omega_\bullet^c)
\leq \P_s(\Omega_\star^c)+ \P_s(\Omega_\sigma^c)+ \P_s(\Omega_M^c)
\leq  C(\theta_5,d)/(n^2\log(n)).$$
From Pythagoras' Equality and the inclusion of $S_\hatmrho$ into $S_\starmr$, we deduce 
$$\|s-\tilde s\|^2_2 = \|s-\hat s_\hatmrho\|_2^2 + \chi^2\hatmrho \leq \|s\|_2^2 + \chi^2\starmr.$$
Therefore, it follows from Cauchy-Scwharz Inequality that 
$$\E_s\left[\|s-\tilde s\|^2_2 \BBone_{\Omega_\bullet}\right] \leq p_\bullet \|s\|_2^2 + \sqrt{p_\bullet\E_s\left[\chi^4\starmr\right]}.$$
Let $\mathscr S_\star$ be some countable and dense subset of $\{t\in S_\starmr \text{ s.t. } \|t\|_2=1\}$. Since $\chi\starmr=\sup_{t\in\mathscr S_\star} |\nu(t)|$, we deduce from Theorem 12 in~\cite{BBLM} that 
$$\sqrt{\E_s\left[\chi^4\starmr\right]}\leq C\left(\E_s\left[\chi^2\starmr\right] + \sigma^2/n + M/n^2\right),$$
where $M$ is any upper-bound for $\sup_{t\in\mathscr S_\star}\max_{1\leq i\leq n} |t(Y_i)-\langle t,s\rangle|$ and $\sigma^2$, any upper-bound for $n\sup_{t\in\mathscr S_\star}\Var_s(t(Y_1)).$ Therefore, we obtain 
$$\sqrt{\E_s\left[\chi^4\starmr\right]}\leq C\left(\frac{\pi(\bb r_\star)|\Lambda(\bb r_\star)| \|s\|_\infty }{\log(n)}+ \frac{\|s\|_\infty}{n} + \frac{\pi(\bb r_\star)|\Lambda(\bb r_\star)|}{n\log(n)}\right) \leq C\frac{\pi(\bb r_\star)|\Lambda(\bb r_\star)|\|s\|_\infty}{\log(n)},$$
hence
\begin{equation}
\E_s\left[\|s-\tilde s\|^2_2 \BBone_{\Omega_\bullet}\right] 
\leq C(\theta_5,d) \frac{\pi(\bb r_\star)|\Lambda(\bb r_\star)|\|s\|^2_\infty}{n\log^{3/2}(n)}.
\label{eq:riskomegacomp}
\end{equation}

Since $\|s\|_\infty\geq 1$, we conclude thanks to~\eqref{eq:riskomega} and~\eqref{eq:riskomegacomp} 
\begin{align}
\E_s\left[\|s-\tilde s\|_2^2\right] 
&\leq C''_1 \|s-s_\mrho\|_2^2 
+C''_2\E_s\left[\|\hat s_\mrho-s_\mrho\|_2^2\right]
+C''_3\pi(\bb r_\star)\frac{D_\mrho}{n}\notag\\
&+\|s\|_\infty \pi(\bb r_\star)|\Lambda(\bb r_\star)| \left( C''_4 \frac{w_\mrho}{n}+C''_5\frac{\Sigma}{n}+C''_6 \frac{\|s\|_\infty}{n\log^{3/2}(n)}\right)
\label{eq:risktotal}
\end{align}
where 
\begin{align*}
C''_1 &= C_1/C_0,\quad C''_2=(1+\theta_6)C_2/C_0,\quad C''_3= 8(4+\theta_6^{-1})C_2/C_0,\\
C''_4 &= (C'_3+C'_4)/C_0,\quad C''_5 =2(C_6+C_7)/C_0,\quad C''_6=C(\theta_5,d).
\end{align*}
Choosing, for all $\mrho\in\Mdeg$, $w_\mrho=L_\mrho |m|$, and taking in~\eqref{eq:risktotal} the minimum over $\mrho\in\Mdeg$ allows to complete the proof.

\section{Proof of Theorem~\ref{adaptselect}}\label{sec:proofadapt}
Let us fix $\bb \sigma, p,p',R,L$ satisfying the assumptions of the theorem and $s\in\mathcal P(\bb \sigma, p,p',R,L)$. For $J=J_\star$, all the partitions given by Theorem~\ref{ApproxAnisotropicC} belong to $\M$, so according to Proposition~\ref{choiceL} and Theorem~\ref{ApproxAnisotropicC} applied with $\bb r=\lfloor {\bb\sigma}\rfloor+1
$,
\begin{align*}
&\E_s\left[\|s-\tilde s\|_2^2\right]\\ 
&\leq C(\kappa'',d,\bb \sigma,p,L) \left(\inf_{k\in\N} \left\{R^22^{-2kH(\bb \sigma)} + w(\bb r_\star) \frac{2^{kd}}{n}\right\} + R^2 2^{-2J_\star d\left(H(\bb \sigma)/d-(1/p-1/2)_+\right) \underline{\bb \sigma}/H(\bb \sigma)}\right).
\end{align*}
In order to minimize approximately the above infimum, we choose 
$$k_\star=\max \{k\in \N \text{ s.t. } w(\bb r_\star) 2^{kd}/n \leq R^22^{-2kH(\bb \sigma)}\}$$
which is well defined since $R^2n/ w(\bb r_\star) \leq 1$, and thus obtain
\begin{align*}
&\E_s\left[\|s-\tilde s\|_2^2\right]\\
&\leq C(\kappa'',d,\bb \sigma,p,L) \left(\left(R\left(n/w(\bb r_\star\right))^{-H(\bb \sigma)/d})\right)^{2d/(d+2H(\bb \sigma))} + R^2 2^{-2J_\star d\left(H(\bb \sigma)/d-(1/p-1/2)_+\right) \underline{\bb \sigma}/H(\bb \sigma)}\right).
\end{align*}
Given the assumptions on $J_\star$ and $R$, the leading term in the right-hand sand is the first one. We then conclude thanks to Propostion~\ref{minimaxonemodel}.

\bibliographystyle{alpha}
\bibliography{BibliographieApprox}

\end{document}